\documentclass[11pt]{article}
\usepackage[T1]{fontenc}
\usepackage[utf8]{inputenc}
\usepackage{lmodern}
\usepackage{amsmath,amssymb,amsthm,mathrsfs}
\usepackage{mathtools}
\usepackage[margin=1in]{geometry}
\usepackage{microtype}
\usepackage{parskip}
\usepackage{needspace}
\usepackage{xcolor}
\usepackage{enumitem}
\usepackage{float}
\usepackage[hidelinks,bookmarksnumbered=true]{hyperref}

\usepackage{tikz}
\usetikzlibrary{arrows.meta}
\usepackage{bookmark}
\newtheorem{theorem}{Theorem}[section]
\newtheorem{lemma}[theorem]{Lemma}
\newtheorem{proposition}[theorem]{Proposition}
\newtheorem*{annulusintro}{Proposition~\ref{prop:annulus}}
\newtheorem{corollary}[theorem]{Corollary}
\newtheorem{standard}[theorem]{Theorem}
\theoremstyle{definition}

\theoremstyle{remark}
\newtheorem{remark}[theorem]{Remark}
\newcommand{\R}{\mathbb{R}}

\setlist{leftmargin=*,itemsep=3pt,topsep=5pt}
\numberwithin{equation}{section}
\AddToHook{env/theorem/before}{\Needspace{4\baselineskip}}
\AddToHook{env/lemma/before}{\Needspace{4\baselineskip}}
\AddToHook{env/proposition/before}{\Needspace{4\baselineskip}}
\AddToHook{env/corollary/before}{\Needspace{4\baselineskip}}
\AddToHook{env/standard/before}{\Needspace{4\baselineskip}}

\hypersetup{pdftitle={On the distance problem of Erd\H{o}s and Davies},pdfauthor={Thang Pham and Chun-Yen Shen},pdfsubject={Distinct distances under a local restriction}}
\title{On the distance problem of Erd\H{o}s and Davies}
\author{Thang Pham \and Chun-Yen Shen}
\date{}

\definecolor{revisionbrown}{RGB}{0,0,255}
\definecolor{revisionpurple}{RGB}{128,0,128}
\definecolor{richrevision}{HTML}{006B45}

\begin{document}
\maketitle
\begin{abstract}

Let $P$ be an $n$-point set in the plane in which every three points
determine three distinct distances. For $p\in P$, let $A_p$ be the set
of squared distances from $p$ to the other points of $P$. We prove that, for every $p\in P$, a positive proportion of
the other points $q\in P$ satisfy $|A_p\setminus A_q|\ge\eta n$,
where $\eta>0$ is absolute and $n$ is sufficiently large. Consequently,
$P$ determines at least $(1+c)n$ distinct distances for an absolute
constant $c>0$. To the best of our knowledge, this is the first fixed
positive proportional improvement over the elementary
lower bound for $D(n,3,3)$.

\par
\end{abstract}
\tableofcontents
\section{Introduction}\label{sec:notation}
{
Erd\H{o}s~\cite[p.~136]{Erdos1973} asked in 1973 whether every set of
$n$ points in the plane containing no isosceles triangle must determine
a superlinear number of distinct distances. The question, which he
attributed to R.~O.~Davies and himself, is listed as Problem~657 in
the Erd\H{o}s problems collection~\cite{Bloom657}. Here degenerate
isosceles triangles are also excluded, so every three points must
determine three distinct distances. Let $D(n,3,3)$ denote the minimum
number of distinct distances determined by such an $n$-point set.
The distances from any fixed point to the remaining $n-1$ points are
all distinct, giving the elementary bound $D(n,3,3)\ge n-1$.
Erd\H{o}s's conjecture~\cite[p.~149]{Erdos1986} is that
\[
 \frac{D(n,3,3)}{n}\longrightarrow\infty
 \qquad\text{as }n\to\infty.
\]

This question has been extended to the more general problem of
estimating $D(n,p,q)$, the minimum number of distinct distances
determined by an $n$-point set in the plane in which every $p$ points
determine at least $q$ distinct distances, where $p\ge3$ and
$1\le q\le\binom{p}{2}$. Fox, Pach, and Suk~\cite{FoxPachSuk}
proved that, for every fixed $p\ge6$,
\[
 D\!\left(n,p,\binom{p}{2}-p+6\right)\ge n^{8/7-o(1)}.
\]
Further lower bounds were obtained by
Pohoata and Sheffer~\cite{PohoataSheffer} and by Fish, Pohoata, and
Sheffer~\cite{FishPohoataSheffer}, through connections with edge
colorings of complete graphs. Complementary constructions were
developed by Fish, Lund, and Sheffer~\cite{FishLundSheffer} using sets of
real numbers whose small subsets determine many distinct differences.
These results give bounds in several parameter ranges, but do not
yield a fixed positive proportional improvement over the elementary
lower bound for $D(n,3,3)$.

In this paper, we prove that there exists an absolute constant $c>0$
such that
\[
 D(n,3,3)\ge(1+c)n
\]
for all sufficiently large $n$. This improves the elementary bound
by a fixed positive proportion of $n$. It follows from the main result
(Theorem~\ref{th:gap}), which compares the pinned distance
set at each point with those at all other points.

Erd\H{o}s~\cite[p.~149]{Erdos1986} also observed that Behrend's
construction of sets of integers without three-term arithmetic
progressions~\cite{Behrend} gives
\[
 D(n,3,3)\le n\exp\!\bigl(C\sqrt{\log n}\bigr)
\]
for an absolute constant $C>0$ and all sufficiently large $n$.
Indeed, on a line, an isosceles triple is precisely a three-term
arithmetic progression. 
For collinear sets satisfying the three-point condition,
Dumitrescu~\cite[Theorem~1]{Dumitrescu} established a
superlinear lower bound of the form $cn(\log n)^\kappa$,
for absolute constants $c,\kappa>0$.
A stronger bound,
$
D(P)\gg n\exp\bigl(c(\log n)^{1/9}\bigr),
$
follows from Ruzsa's difference-set inequality \cite{RUUU} and the
estimate of Bloom and Sisask for sets without three-term
arithmetic progressions~\cite{BloomSisask}.

\par}

Before stating our main result, we first introduce some notations first. A finite set $P\subset\R^2$ is \emph{proper} if
\[
 |p-x|\ne |p-y|\qquad(p,x,y\in P\text{ pairwise distinct}).
\]

Thus, properness is equivalent to this three-point condition.
For a nonempty proper set $P$, let $n=|P|$, and let $D(P)$ be the number
of distinct positive distances determined by $P$. Note that squaring the distances
does not change this number.
The \emph{pinned squared-distance set} at $p\in P$ is
\[
 A_p=\{|p-x|^2:x\in P\setminus\{p\}\}.
\]

Properness gives $|A_p|=|P|-1$. For $p,q\in P$, define the \emph{pinned disagreement}
and its average by
\[
 d_P(p,q)=|A_p\setminus A_q|,
 \qquad
 \overline d_P(p)=\frac1n\sum_{q\in P}d_P(p,q),
\]
\[
 d_{\mathrm{av}}(P)=\frac1{n^2}\sum_{p,q\in P}d_P(p,q),
 \qquad d_{\max}(P)=\max_{p,q\in P}d_P(p,q).
\]
Since $|A_p|=|A_q|$, one has
$d_P(p,q)=\tfrac12|A_p\mathbin{\triangle}A_q|$, where
$\triangle$ denotes \emph{symmetric difference}. Thus, $d_P$ is
symmetric and satisfies the triangle inequality. We omit the subscript $P$ when it is clear.
The averages $\overline d_P(p)$ and $d_{\mathrm{av}}(P)$ are used only for nonempty sets.
If $S\subset P$, the symbols $A_{p,S}$ and $d_S$ use both
pins and targets in $S$.

For $t>0$, define the \emph{ordered multiplicity} and
\emph{positive distance energy} by
\[
 \nu_t(P)=|\{(x,y)\in P^2:|x-y|^2=t\}|,
 \qquad E(P)=\sum_{t>0}\nu_t(P)^2.
\]

Since $\sum_{t>0}\nu_t(P)=n(n-1)$, the Cauchy--Schwarz
inequality gives
\begin{equation}\label{eq:cs-distances}
 D(P)\ge\frac{n^2(n-1)^2}{E(P)}
 \qquad(n\ge2).
\end{equation}
\par

Every positive distance class is a \emph{matching}: its edges have
no common endpoint. Hence, $\nu_t(P)$ equals the number of pins whose
pinned distance set contains $t$. An \emph{equal-distance quadruple} is an
ordered quadruple
$(p_1,p_2,q_1,q_2)$ of points in $\R^2$ satisfying
$|p_1-p_2|=|q_1-q_2|>0$. The ordered pairs $(p_1,p_2)$ and
$(q_1,q_2)$ are its \emph{source segment} and \emph{target segment},
respectively. An isometry $f$ \emph{realizes} the quadruple if
$f(p_i)=q_i$ for $i=1,2$. Points from different segments may coincide,
but the two endpoints within each segment are distinct. Simultaneously
reversing both segments gives a different ordered quadruple.
Thus, $E(P)$ counts the equal-distance quadruples in $P^4$.
For finite sets $P,T$, the \emph{mixed distance energy} is
\[
 E(P,T)=|\{(p,q,x,y)\in P^2\times T^2:|p-x|^2=|q-y|^2\}|.
\]
This definition includes zero distances; in particular,
$E(P,P)=E(P)+|P|^2$.

For quantities $f$ and $g>0$, we write $f=O(g)$ or $f\ll g$ if
$|f|\le Cg$ for some constant $C>0$ in the range under consideration.
We write $f=o(g)$ if $f/g\to0$, and $f\sim g$ if $f/g\to1$,
in the indicated limit. Constants implicit in $O$ and $\ll$ are
absolute unless a dependence is indicated. Subscripts in $O$, $\ll$,
and $o$ list fixed parameters on which the implicit constant or
convergence threshold may depend. All little-oh statements about
families of configurations are uniform unless a fixed sequence is
explicitly specified.
\par

Our main theorem, which studies the interaction between pinned distance sets
at different points, reads as follows.

\begin{theorem}\label{th:gap}
There are absolute constants $\eta>0$ and $n_0$ such that every proper
set $P\subset\R^2$ with $n=|P|\ge n_0$ satisfies, for every $p\in P$,
\[
 \overline d_P(p)\ge\eta n,\qquad
 D(P)\ge n-1+\eta n,\qquad
 E(P)\le (1-\eta)n^3.
\]
\end{theorem}

\begin{corollary}\label{cor:distinct-distances}
There are absolute constants $c>0$ and $n_1$ such that
\[
 D(n,3,3)\ge(1+c)n\qquad(n\ge n_1).
\]
\end{corollary}

\begin{corollary}\label{cor:many-pins}
There are absolute constants $\eta>0$ and $n_0$ such that every proper
set $P\subset\R^2$ with $n=|P|\ge n_0$ satisfies, for every $p\in P$,
\[
 \bigl|\{q\in P\setminus\{p\}:|A_p\setminus A_q|\ge\eta n\}\bigr|
 \ge\eta n.
\]
\end{corollary}

The bound for each pin in Theorem~\ref{th:gap} has optimal linear
order, since $\overline d_P(p)\le(n-1)(n-2)/n$.
Indeed, for $q\ne p$, the common squared distance $|p-q|^2$ belongs
to both $A_p$ and $A_q$, so $|A_p\setminus A_q|\le n-2$.

The proof of Theorem~\ref{th:gap} has three main ingredients,
Propositions~\ref{prop:energy-pruning}, \ref{prop:polygon-matching},
and~\ref{prop:annulus}. The central geometric ingredient,
Proposition~\ref{prop:annulus}, which is of independent interest, gives subcubic energy for proper
sets in shrinking annuli under \emph{angular nonconcentration}:
short angular intervals contain uniformly small proportions of the
points. No relation is imposed between the rate at which the annuli
shrink and the number of points. It is worth noting that geometric convergence to a circle
alone does not yield this energy estimate.

The proof also uses a mixed distance-energy estimate for
points on conics. Pach and de Zeeuw proved a superlinear
distinct-distance bound for point sets supported on
noncircular irreducible conics. Our Theorem~\ref{thm:conic-target-energy}
instead gives
\[
E(P,T)=o(|P|^3)
\]
when $P$ is a proper planar set and $T\subset P$ lies on
an irreducible real conic, including a circle.
The mixed-energy form, which allows the pins to range
over the full set $P$, enables us to remove concentrations
on conics in the proof of Proposition~\ref{prop:annulus}.

We state Proposition~\ref{prop:annulus} here. For $x\ne0$, let
$\arg x\in\R/(2\pi\mathbb Z)$ denote its angle. An \emph{angular interval}
is an arc in this angle circle, and $|I|$ denotes its length in radians,
between zero and $2\pi$.

\Needspace{8\baselineskip}
\begin{annulusintro}
Let $P_j\subset\R^2$ be proper sets with $|P_j|\to\infty$.
Suppose that $0<h_j<1/2$, $h_j\to0$, and $P'_j\subset P_j$ satisfy
\[
 |P_j\setminus P'_j|=o(|P_j|),\qquad
 \bigl||x|-1\bigr|\le h_j\qquad(x\in P'_j).
\]
Suppose also that, for every $\varepsilon>0$, there are
$0<t\le2\pi$ and an integer $j_0$ such that, for every $j\ge j_0$
and every angular interval $I$ with $|I|\le t$,
\[
 \bigl|\{x\in P'_j:\arg x\in I\}\bigr|\le\varepsilon|P_j|.
\]

Then
\[
 E(P_j)=o(|P_j|^3),\qquad
 \frac{D(P_j)}{|P_j|}\longrightarrow\infty.
\]
\end{annulusintro}

Neither properness nor angular nonconcentration can be omitted
from Proposition~\ref{prop:annulus}.
Examples illustrating this will be given in
Section~\ref{sec:annulus-examples}.

\paragraph{Organization.}
Section~\ref{sec:main-proof} deduces Theorem~\ref{th:gap} and its
corollaries from three propositions. Section~\ref{sec:proposition-proofs}
opens with the proof of the energy and pruning estimates in
Proposition~\ref{prop:energy-pruning}. Subsection~\ref{sec:polygon-matching}
proves the approximation by a regular polygon.
Subsections~\ref{tri:upd:sec:rich} and~\ref{sec:curve-targets} bound
the number of isometries with many matches and the mixed distance energy
for targets on a conic, respectively. Subsection~\ref{sec:annulus-proof}
uses these bounds to prove Proposition~\ref{prop:annulus} on shrinking
annuli. Section~\ref{sec:annulus-examples} gives examples showing that
both properness and angular nonconcentration are necessary for its
energy bound. Section~\ref{sec:polynomial-conjecture} asks which
symmetric polynomials can take only linearly many values on arbitrarily
large admissible sets.
\par

\section{Proof of Theorem~\ref{th:gap} from three propositions}\label{sec:main-proof}

We deduce the theorem from three propositions. The first relates
distance energy to the differences between pinned distance sets and
allows us to discard a small exceptional set. The second matches the
remaining points to a regular polygon. The third gives subcubic
energy from the resulting radial and angular bounds.

In proving Theorem~\ref{th:gap}, we may assume that $P$ has
centroid zero and diameter two. Indeed, for $|P|\ge2$, let
$b=|P|^{-1}\sum_{p\in P}p$. The similarity
\[
 x\longmapsto \frac{2(x-b)}{\operatorname{diam}P}
\]
sends $P$ to a set with these properties. It multiplies all
squared distances by the same positive factor and therefore
preserves properness, $D(P)$, $E(P)$, and all pinned disagreement
quantities.
\par

The first ingredient shows that if the pinned distance set at one point
differs from those at the other points in only $o(|P|)$ elements on average,
we can delete $o(|P|)$ points while retaining that point. In the remaining
set $S$, every pair of pinned distance sets, now computed within $S$,
differs in $o(|S|)$ elements, and consequently
$E(S)=(1-o(1))|S|^3$.
\par

\begin{proposition}\label{prop:energy-pruning}
\begin{enumerate}[label=\textup{(\roman*)}]
\item Let $P\subset\R^2$ be a nonempty proper set and let $n=|P|$.
For every $p\in P$,
\[
 E(P)=n^2(n-1-d_{\mathrm{av}}(P)),\qquad
 d_{\mathrm{av}}(P)\le2\overline d_P(p).
\]
Consequently,
\[
 n^2(n-1-2\overline d_P(p))\le E(P)\le n^2(n-1),
 \qquad \overline d_P(p)\le D(P)-(n-1).
\]

\item Let $P_j\subset\R^2$ be proper sets with $|P_j|\to\infty$,
and suppose that $p_j\in P_j$ satisfy
$\overline d_{P_j}(p_j)=o(|P_j|)$. There are subsets
$S_j\subset P_j$ containing $p_j$ such that
\[
 |S_j|\sim|P_j|,\qquad d_{\max}(S_j)=o(|S_j|),\qquad
 E(S_j)=(1-o(1))|S_j|^3.
\]
\end{enumerate}
\end{proposition}

For the next two proposition, let
\[
 S^1=\{x\in\R^2:|x|=1\}.
\]

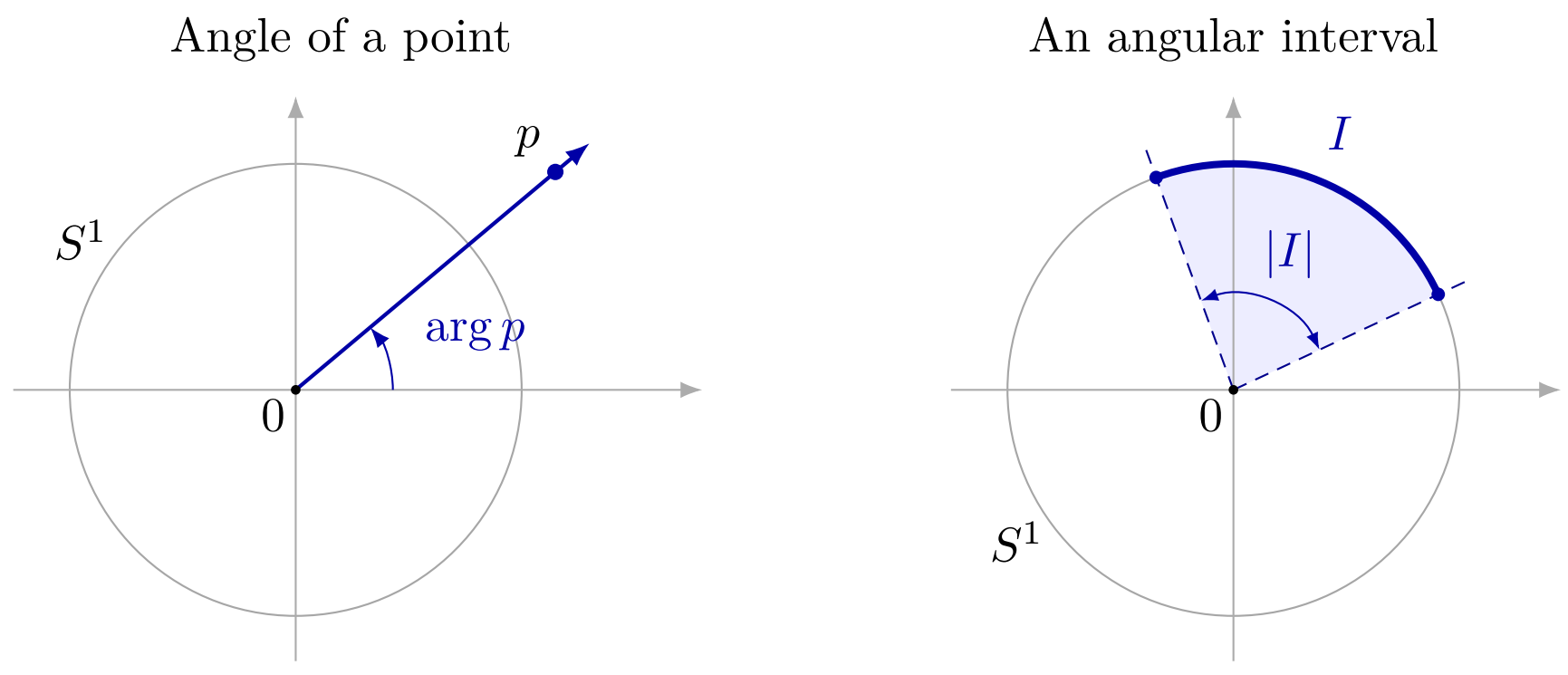
\begin{figure}[H]
\centering
\begin{tikzpicture}[
  x=1.65cm,y=1.65cm,
  font=\small,
  axis/.style={-{Latex[length=1.6mm]},gray!65,thin},
  circleline/.style={gray!70,thin},
  ray/.style={blue!65!black,thick}
]
\begin{scope}
  \node at (.2,1.55) {Angle of a point};
  \draw[circleline] (0,0) circle[radius=1];
  \draw[axis] (-1.25,0) -- (1.8,0);
  \draw[axis] (0,-1.2) -- (0,1.3);
  \draw[ray,-{Latex[length=1.8mm]}] (0,0) -- (40:1.7);
  \fill[blue!65!black] (40:1.5) circle[radius=1.7pt];
  \node[above left,inner sep=3pt] at (40:1.5) {$p$};
  \draw[blue!65!black,-{Latex[length=1.5mm]},thin]
    (.43,0) arc[start angle=0,end angle=40,radius=.43];
  \node[blue!65!black,inner sep=1pt] at (17:.83) {$\arg p$};
  \fill (0,0) circle[radius=1pt];
  \node[below left,inner sep=2pt] at (0,0) {$0$};
  \node at (145:1.16) {$S^1$};
\end{scope}
\begin{scope}[shift={(4.15,0)}]
  \node at (0,1.55) {An angular interval};
  \fill[blue!7] (0,0) -- (25:1)
    arc[start angle=25,end angle=110,radius=1] -- cycle;
  \draw[circleline] (0,0) circle[radius=1];
  \draw[axis] (-1.25,0) -- (1.45,0);
  \draw[axis] (0,-1.2) -- (0,1.3);
  \draw[blue!55!black,densely dashed,thin]
    (0,0) -- (25:1.17) (0,0) -- (110:1.17);
  \draw[blue!65!black,line width=1.5pt]
    (25:1) arc[start angle=25,end angle=110,radius=1];
  \fill[blue!65!black] (25:1) circle[radius=1.4pt];
  \fill[blue!65!black] (110:1) circle[radius=1.4pt];
  \draw[blue!65!black,
    {Latex[length=1.25mm]}-{Latex[length=1.25mm]},thin]
    (25:.42) arc[start angle=25,end angle=110,radius=.42];
  \node[blue!65!black,fill=blue!7,inner sep=1pt] at (67.5:.65) {$|I|$};
  \node[blue!65!black] at (67.5:1.23) {$I$};
  \fill (0,0) circle[radius=1pt];
  \node[below left,inner sep=2pt] at (0,0) {$0$};
  \node at (215:1.17) {$S^1$};
\end{scope}
\end{tikzpicture}
\caption{The angle $\arg p$ is measured counterclockwise from the
positive horizontal axis, modulo $2\pi$. On the right, the highlighted
arc on $S^1$ represents an angular interval $I$ in
$\R/(2\pi\mathbb Z)$, and the central angle is its length $|I|$ in radians.
The shaded sector indicates the corresponding directions.}
\label{fig:angle-angular-interval}
\end{figure}

The second ingredient shows that, under the normalization above,
$d_{\max}(P)=o(|P|)$ implies that $P$ can be matched bijectively
to the vertices of a regular $|P|$-gon on $S^1$, with maximum
displacement $o(1)$. Consequently, all points lie in a shrinking
annulus about $S^1$, and their angular distribution converges
to the uniform distribution, uniformly over angular intervals.
\par

\begin{proposition}\label{prop:polygon-matching}
There is an absolute constant $C\ge1$ with the following property.
Let $P\subset\R^2$ be proper, with $|P|\ge3$, centroid zero, and
diameter two. Set
\[
 \rho=C\frac{d_{\max}(P)+1}{|P|}.
\]
There is a bijection $p\mapsto v(p)$ from $P$ to the vertices of a
regular $|P|$-gon on $S^1$ such that
\[
 \max_{p\in P}|p-v(p)|\le\rho.
\]

In particular,
\[
 \bigl||p|-1\bigr|\le\rho\qquad(p\in P).
\]
If $\rho\le1/4$, then every angular interval $I$ satisfies
\[
 \left|\frac{\bigl|\{p\in P:\arg p\in I\}\bigr|}{|P|}
       -\frac{|I|}{2\pi}\right|
 \le C\rho+\frac{2}{|P|}
 \le\left(C+\frac{2}{C}\right)\rho.
\]
\end{proposition}

The matching follows from the finite circular stability theorem of
Pham, Tran, and Shen~\cite[Corollary~1.3]{PhamTShen}. Its application to
proper sets and the angular estimate are proved in
Subsection~\ref{sec:polygon-matching}.

\begin{proposition}\label{prop:annulus}
Let $P_j\subset\R^2$ be proper sets with $|P_j|\to\infty$.
Suppose that $0<h_j<1/2$, $h_j\to0$, and $P'_j\subset P_j$ satisfy
\[
 |P_j\setminus P'_j|=o(|P_j|),\qquad
 \bigl||x|-1\bigr|\le h_j\qquad(x\in P'_j).
\]

Suppose also that, for every $\varepsilon>0$, there are
$0<t\le2\pi$ and an integer $j_0$ such that, for every $j\ge j_0$
and every angular interval $I$ with $|I|\le t$,
\begin{equation}\label{eq:angular-nonconcentration}
 \bigl|\{x\in P'_j:\arg x\in I\}\bigr|\le\varepsilon|P_j|.
\end{equation}

Then
\[
 E(P_j)=o(|P_j|^3),\qquad
 \frac{D(P_j)}{|P_j|}\longrightarrow\infty.
\]
\end{proposition}

\begin{proof}[Proof of Theorem~\ref{th:gap}]
It suffices to prove the asserted uniform lower bound for
$\overline d_P(p)$. Indeed, suppose that this bound holds with
$0<\eta<1$ for every proper set $P$ with $n=|P|\ge n_0$ and every
$p\in P$. Proposition~\ref{prop:energy-pruning}\textup{(i)} then gives
\[
 D(P)\ge n-1+\overline d_P(p)\ge n-1+\eta n,
 \qquad
 d_{\mathrm{av}}(P)=\frac1n\sum_{p\in P}\overline d_P(p)\ge\eta n.
\]
The energy identity in the same proposition yields
\[
 E(P)\le n^2(n-1-\eta n)\le(1-\eta)n^3.
\]
These inequalities give the remaining conclusions of the theorem.

We now prove the uniform lower bound. Suppose, for a contradiction, that no such constants $\eta>0$
and $n_0$ exist. Then, for every $\eta>0$ and every integer
$n_0\ge2$, there are a proper set $P$ with $|P|\ge n_0$
and a point $p\in P$ such that
$\overline d_P(p)<\eta|P|$.
Taking $\eta=1/j$ and $n_0=j$, we obtain, for each integer
$j\ge2$, a proper set $P_j$ and a point $p_j\in P_j$ satisfying
\[
 |P_j|\ge j,\qquad
 \overline d_{P_j}(p_j)<\frac{|P_j|}{j}.
\]
Proposition~\ref{prop:energy-pruning}\textup{(ii)} gives subsets
$S_j\subset P_j$ such that
\[
 |S_j|\sim|P_j|,\qquad d_{\max}(S_j)=o(|S_j|),\qquad
 E(S_j)=(1-o(1))|S_j|^3.
\]
We will use the other two propositions to show that these same sets
satisfy $E(S_j)=o(|S_j|^3)$.

For all sufficiently large $j$, normalize each $S_j$ as above and
retain the same notation; the preceding disagreement and energy
estimates are preserved.
\par

Let $C$ be the constant in Proposition~\ref{prop:polygon-matching}
and set
\[
 \rho_j=C\frac{d_{\max}(S_j)+1}{|S_j|}\longrightarrow0.
\]
For all sufficiently large $j$, we have $|S_j|\ge3$ and
$0<\rho_j\le1/4$. That proposition therefore gives
\[
 \bigl||p|-1\bigr|\le\rho_j\qquad(p\in S_j)
\]

and, for every angular interval $I$,
\[
 \frac{1}{|S_j|}\bigl|\{p\in S_j:\arg p\in I\}\bigr|
 \le\frac{|I|}{2\pi}+C\rho_j+\frac{2}{|S_j|}.
\]

To verify the angular hypothesis of Proposition~\ref{prop:annulus},
fix $\varepsilon>0$ and choose $0<t\le2\pi$ with
$t/(2\pi)<\varepsilon/2$. For all sufficiently large $j$,
$C\rho_j+2/|S_j|<\varepsilon/2$. The displayed estimate then gives
\eqref{eq:angular-nonconcentration} for every interval with $|I|\le t$.
The radial hypothesis holds with $h_j=\rho_j<1/2$, and no points
need to be discarded.

Proposition~\ref{prop:annulus}, applied with
$P_j=P'_j=S_j$ and $h_j=\rho_j$, now gives $E(S_j)=o(|S_j|^3)$.
This contradicts the energy estimate from pruning: the ratio
$E(S_j)/|S_j|^3$ cannot tend both to zero and to one.
Thus, the required uniform lower bound holds for some $\eta>0$ and
$n_0$. Decreasing $\eta$ to be less than one and applying the initial
reduction proves all conclusions of Theorem~\ref{th:gap}.
\end{proof}

\begin{proof}[Proof of Corollary~\ref{cor:distinct-distances}]
Let $\eta$ and $n_0$ be as in Theorem~\ref{th:gap}. Set $c=\eta/2$
and choose an integer $n_1\ge\max\{n_0,2/\eta\}$. Every proper
set $P$ of $n\ge n_1$ points satisfies
\[
 D(P)\ge n-1+\eta n\ge(1+c)n.
\]
Taking the minimum over all such sets proves the corollary.
\end{proof}

\begin{proof}[Proof of Corollary~\ref{cor:many-pins}]
Let $\eta_0>0$ and $n_0$ be constants supplied by Theorem~\ref{th:gap},
and set $\eta=\eta_0/2$. Fix a proper set $P$ with $n=|P|\ge n_0$
and a point $p\in P$, and let
\[
 Q_p=\{q\in P\setminus\{p\}:|A_p\setminus A_q|\ge\eta n\}.
\]
Since $|A_p\setminus A_q|\le n-1$ for every $q\in P$, the theorem gives
\[
 2\eta n^2
 \le\sum_{q\in P}|A_p\setminus A_q|
 \le |Q_p|n+\eta n^2.
\]
Thus, $|Q_p|\ge\eta n$, which proves the corollary.
\end{proof}

\section{Proofs of the three propositions}\label{sec:proposition-proofs}
We first establish the energy estimates and the approximation by a
regular polygon. We then prove the estimates for Euclidean motions
and conics that are needed for the annulus argument.
The Proposition~\ref{prop:energy-pruning} is proved by counting intersections
of pinned distance sets and then deleting the points whose pinned
distance sets differ too much from that at a fixed point. The following
restriction estimates control
the effect of those deletions. The energy estimate will also be used
in Subsection~\ref{sec:annulus-proof}.

\begin{lemma}\label{lem:restriction}
Let $P$ be proper and $S\subset P$. For $x,y\in S$,
\[
 d_S(x,y)\le d_P(x,y)+|P\setminus S|.
\]

Moreover,
\begin{equation}\label{eq:deletion-energy}
 0\le E(P)-E(S)\le4|P\setminus S|\,|P|^2.
\end{equation}
\end{lemma}

\begin{proof}
For $x,y\in S$, the inclusions $A_{x,S}\subset A_x$ and
$A_{y,S}\subset A_y$ give
\[
 A_{x,S}\setminus A_{y,S}
 \subset (A_x\setminus A_y)
          \cup(A_y\setminus A_{y,S}).
\]

Properness at $y$ gives
$|A_y\setminus A_{y,S}|=|P\setminus S|$. Consequently,
\[
 d_S(x,y)
 \le d_P(x,y)+|A_y\setminus A_{y,S}|
 =d_P(x,y)+|P\setminus S|.
\]

The difference $E(P)-E(S)$ counts the equal-distance
quadruples with at least one endpoint in $P\setminus S$, so it is
nonnegative. Fix a deleted point in one of the four endpoint positions.
There are at most $|P|^2$ choices for the other endpoint of its segment
and the first endpoint of the other segment. Properness determines
the last endpoint in at most one way, because the common distance
is positive. Summing over the four positions and the deleted points gives
\[
 E(P)-E(S)\le4|P\setminus S|\,|P|^2.
\]

This proves~\eqref{eq:deletion-energy} and
completes the proof of the lemma.
\end{proof}

\begin{proof}[Proof of Proposition~\ref{prop:energy-pruning}]
We first prove part~(i). Let $P$ be a nonempty proper set and write
$n=|P|$. For each positive squared distance $t$, properness gives at
most one point at that squared distance from each pin. Thus,
$\nu_t(P)$ equals the number of pinned distance sets $A_x$ containing $t$.
Counting the pairs of pinned distance sets containing each $t$ gives
\[
 \begin{aligned}
 E(P)&=\sum_{x,y\in P}|A_x\cap A_y|\\
 &=n^2(n-1)-\sum_{x,y\in P}d_P(x,y)
 =n^2(n-1-d_{\mathrm{av}}(P)).
 \end{aligned}
\]
Here $|A_x\cap A_y|=n-1-d_P(x,y)$ because every pinned
distance set has $n-1$ elements. For a fixed $p\in P$, the triangle inequality and
symmetry of $d_P$ imply
\[
 d_{\mathrm{av}}(P)
 \le\frac1{n^2}\sum_{x,y\in P}
       \bigl(d_P(x,p)+d_P(p,y)\bigr)
 =2\overline d_P(p).
\]
Together with $d_{\mathrm{av}}(P)\ge0$, this proves both energy bounds.

The union of any two pinned distance sets has at most $D(P)$ elements. Hence, for
every $x\in P$,
\[
 |A_p\setminus A_x|=|A_p\cup A_x|-|A_x|
 \le D(P)-(n-1).
\]
Averaging over $x$ gives the claimed bound for
$\overline d_P(p)$.
This completes part~(i).

For part~(ii), we construct a suitable subset for every $p\in P$.
For a positive threshold $s$, there are at most
$n\overline d_P(p)/s$ points $x$ with
$d_P(p,x)>s$. For any two remaining points, the triangle
inequality gives $d_P(x,y)\le2s$. We balance these bounds by setting
\[
 s=\sqrt{n(\overline d_P(p)+1)},\qquad
 S=\{x\in P:d_P(p,x)\le s\}.
\]

The added $1$ ensures $s>0$ even when the mean disagreement is zero.
Also, $p\in S$, since $d_P(p,p)=0$. Every deleted pin contributes
more than $s$ to the sum defining the mean, so
\[
 |P\setminus S|s\le\sum_{x\in P\setminus S}d_P(p,x)
 \le n\overline d_P(p),
 \qquad
 |P\setminus S|\le\frac{n\overline d_P(p)}s\le s.
\]

The last inequality uses $n\overline d_P(p)\le s^2$.
For $x,y\in S$, Lemma~\ref{lem:restriction} and the triangle
inequality now give
\[
 d_S(x,y)
 \le d_P(x,p)+d_P(p,y)+|P\setminus S|
 \le2s+|P\setminus S|\le3s.
\]

Thus, $d_{\max}(S)\le3s$.

Along the sequence in part~(ii), apply this construction to each
$(P_j,p_j)$. Suppressing the index, the hypothesis gives
\[
 \frac{s}{n}=\sqrt{\frac{\overline d_P(p)}n+\frac1n}
 \longrightarrow0.
\]
It follows that $|P\setminus S|=o(n)$, $|S|\sim n$, and
$d_{\max}(S)=o(|S|)$.
Since $0\le d_{\mathrm{av}}(S)\le d_{\max}(S)$, the identity proved in part~(i)
yields
\[
 \frac{E(S)}{|S|^3}
 =1-\frac1{|S|}-\frac{d_{\mathrm{av}}(S)}{|S|}\longrightarrow1.
\]
The sets $S_j$ therefore have every property asserted in part~(ii).
This completes the proof of the proposition.
\end{proof}

\subsection{Approximation by a regular polygon}\label{sec:polygon-matching}
We deduce Proposition~\ref{prop:polygon-matching} from the finite
circular stability theorem of Pham, Tran, and
Shen~\cite[Corollary~1.3]{PhamTShen}. We then estimate radial errors
and angular point counts under the resulting matching.

For $p\in\R^2$ and $r\ge0$, let
$\overline{B}(p,r)=\{x\in\R^2:|x-p|\le r\}$.
\par

\Needspace{11\baselineskip}
\begin{standard}\label{lit:polygon-matching}
Let $P\subset\R^2$ have $n\ge3$ points, diameter $L$, and \emph{centroid}
$b=n^{-1}\sum_{p\in P}p$. Suppose that $k\ge0$ and
\[
 \bigl|\,|P\cap\overline{B}(p,r)|-|P\cap\overline{B}(q,r)|\,\bigr|
 \le k\qquad(p,q\in P,\ r\ge0).
\]
There is a regular $n$-gon centered at $b$, of circumradius $L/2$,
and a labeling of its vertices $v_1,\ldots,v_n$ and of the points
$P=\{p_1,\ldots,p_n\}$ such that
\[
 \max_{1\le i\le n}|p_i-v_i|\le CL\,\frac{k+1}{n},
\]
where $C$ is absolute.
\end{standard}

The following elementary estimate applies to any finite set matched
closely to a regular polygon. No properness assumption is needed.

\begin{lemma}\label{lem:polygon-counts}
Let $v_1,\ldots,v_n$ be the vertices of a regular $n$-gon on the
unit circle. Suppose $P=\{p_1,\ldots,p_n\}$ satisfies
$|p_i-v_i|\le\rho\le1/4$ for every $i$. Then
\[
 \bigl||p_i|-1\bigr|\le\rho\qquad(1\le i\le n),
\]
and every angular interval $I$ satisfies
\begin{equation}\label{eq:polygon-counts}
 \bigl|\{p\in P:\arg p\in I\}\bigr|
 \le \frac{n|I|}{2\pi}+Cn\rho+2,
\end{equation}
where $|I|$ is its angular length and $C$ is absolute.
\end{lemma}
\begin{proof}
The reverse triangle inequality gives
$\bigl||p_i|-1\bigr|\le|p_i-v_i|\le\rho$. In particular,
$|p_i|\ge3/4$, so all the angles in the statement are defined.
To compare the angular positions, observe that
\[
 \left|\frac{p_i}{|p_i|}-v_i\right|
 \le\bigl||p_i|-1\bigr|+|p_i-v_i|\le2\rho.
\]
A chord subtending an angle $u\in[0,\pi]$ has length
$2\sin(u/2)\ge2u/\pi$. Hence, the smaller angular distance between
$p_i$ and $v_i$ is at most $\pi\rho$.

If $\arg p_i\in I$, then $\arg v_i$ belongs to the arc obtained by
enlarging $I$ at each end by $\pi\rho$. An arc of length
$t\le2\pi$ contains at most $nt/(2\pi)+2$ vertices of a regular
$n$-gon, including its endpoints. When $|I|+2\pi\rho\le2\pi$, this gives
\[
 \bigl|\{p\in P:\arg p\in I\}\bigr|
 \le\frac{n(|I|+2\pi\rho)}{2\pi}+2
 =\frac{n|I|}{2\pi}+n\rho+2.
\]

If $|I|+2\pi\rho>2\pi$, the bound by $n$ is smaller than the same
right-hand side. This proves~\eqref{eq:polygon-counts} and completes
the proof of the lemma.
\end{proof}

\begin{proof}[Proof of Proposition~\ref{prop:polygon-matching}]
We bound the differences between disk counts in terms of
$d_{\max}(P)$ so that we can apply
Theorem~\ref{lit:polygon-matching}. For every $p\in P$ and $r\ge0$,
properness gives
\[
 |P\cap\overline{B}(p,r)|=1+|A_p\cap[0,r^2]|.
\]
For $p,q\in P$, we have
\[
 \begin{aligned}
 |A_p\cap[0,r^2]|-|A_q\cap[0,r^2]|
 &=|(A_p\setminus A_q)\cap[0,r^2]|\\
 &\quad-|(A_q\setminus A_p)\cap[0,r^2]|.
 \end{aligned}
\]

Each term on the right lies between zero and
$d_{\max}(P)$. Hence,
\[
 \bigl|\,|P\cap\overline{B}(p,r)|-|P\cap\overline{B}(q,r)|\,\bigr|
 \le d_{\max}(P).
\]
Choose the absolute constant $C\ge1$ in the proposition to dominate
both the matching constant in Theorem~\ref{lit:polygon-matching}
with $L=2$ and the angular-count constant in
Lemma~\ref{lem:polygon-counts}. The theorem applies with $k=d_{\max}(P)$.
Since the centroid is zero and the diameter is two, its polygon lies
on $S^1$ and the resulting bijection has displacement at most
$\rho=C(d_{\max}(P)+1)/|P|$.

Since $|v(p)|=1$, the reverse triangle inequality give
\[
 \bigl||p|-1\bigr|\le |p-v(p)|\le\rho\qquad(p\in P).
\]
\par
If $\rho\le1/4$, Lemma~\ref{lem:polygon-counts} applied to this
bijection gives
\[
 \bigl|\{p\in P:\arg p\in I\}\bigr|
 \le\frac{|P|\,|I|}{2\pi}+C|P|\rho+2.
\]

Division by $|P|>0$ gives the upper bound. Applying the same estimate
to a complementary arc gives the corresponding lower bound.
Since $\rho\ge C/|P|$, we obtain the stated discrepancy estimate.
This proves the proposition.
\end{proof}
In the next subsections, we will develop tools to prove Proposition 2.3. 
\subsection{Isometries with many matches}
\label{app:triangle-energy}\label{tri:sec:triangle-background}
\label{tri:upd:sec:rich}

Let $\mathrm{SE}(2)$ denote the group of orientation-preserving
isometries of the plane, also called \emph{direct isometries}. For
finite sets $P,Q$, a \emph{match} of an isometry $g$ is a point
$p\in P$ such that $gp\in Q$. We prove that only $o(n)$ isometries
can have a fixed positive proportion of $n$ matches between two
proper sets of size at most $n$. This estimate will control
concentrations on conics in the proof of Proposition~\ref{prop:annulus}. For related structural results concerning Euclidean motions,
see Mansfield and Passant~\cite{MansfieldPassant}.
Our argument gives a power-saving bound for proper sets
using a finite averaging argument and Roth's theorem.

\begin{theorem}\label{tri:upd:thm:rich}
For every fixed $\varepsilon>0$, there are constants
$c_\varepsilon,K_\varepsilon>0$ such that, for all sufficiently
large $n$, every proper planar set $P$ and every planar set $Q$,
both of size at most $n$, satisfy
\begin{equation}\label{tri:upd:eq:rich-power-saving}
\bigl|\{g\in SE(2):|Q\cap gP|\ge\varepsilon n\}\bigr|
\le K_\varepsilon n^{1-c_\varepsilon}.
\end{equation}
The same bound holds for orientation-reversing isometries.
\end{theorem}

\begin{remark}[Outline of the proof of Theorem~\ref{tri:upd:thm:rich}]
\label{rem:rich-motions-strategy}
The proof counts triples consisting of a motion and two lists
of motions whose averaged indicator functions are close after
applying the first motion. For each fixed motion, a finite
averaging lemma provides many pairs of lists. Conversely,
properness and Roth's theorem show that each fixed pair of
lists works for only a bounded number of motions.
Comparing these two counts proves the theorem.
\end{remark}

For a finitely supported real-valued function $u$ on the plane and
an isometry $g$, write
\[
 \|u\|_n^2=\frac1n\sum_x|u(x)|^2,
 \qquad (gu)(x)=u(g^{-1}x).
\]
Write $\mathbf1_A$ for the indicator of a set $A$.
Every isometry preserves this norm, and $g\mathbf1_A=\mathbf1_{gA}$.

We will use the following elementary estimate. If $P$ is proper,
$|P|\le n$, and $\rho n\ge2$, then
\begin{equation}\label{tri:upd:eq:self-rich}
 \bigl|\{f\in\mathrm{SE}(2):|P\cap fP|\ge\rho n\}\bigr|
 \le \frac{2}{\rho^2}n.
\end{equation}
To see this, count equal-distance quadruples $(p_1,p_2,q_1,q_2)\in P^4$
with $|p_1-p_2|=|q_1-q_2|>0$. There are at most $n^3$: after
$p_1,p_2,q_1$ are chosen, properness determines $q_2$ in at most one
way. Each quadruple determines exactly one direct isometry taking
$p_i$ to $q_i$ for $i=1,2$, including when that isometry is a
translation. A motion $f$ with $r=|P\cap fP|\ge\rho n$ matches accounts for
$r(r-1)\ge\rho^2n^2/2$ quadruples. This proves~\eqref{tri:upd:eq:self-rich}.

\begin{lemma}\label{tri:upd:lem:averages}
Let $M\ge1$ be an integer, and let $u_1,\ldots,u_M$ be vectors of norm at most one in a real inner product
space. If $t>0$ and $\ell\ge4/t^2$ is an integer, then at least half of the
$M^{2\ell}$ ordered pairs of $\ell$-tuples $(i_1,\ldots,i_\ell)$ and
$(j_1,\ldots,j_\ell)$ with entries in $\{1,\ldots,M\}$ satisfy
\[
 \left\|\frac1\ell\sum_{r=1}^{\ell}u_{i_r}
       -\frac1\ell\sum_{r=1}^{\ell}u_{j_r}\right\|\le t.
\]
\end{lemma}
\begin{proof}
Put $\bar u=M^{-1}\sum_{i=1}^M u_i$. Averaging the squared norm over all
ordered pairs of tuples and expanding the inner products gives
\[
 \frac{2}{\ell}\left(\frac1M\sum_{i=1}^M\|u_i\|^2
       -\|\bar u\|^2\right)\le\frac2\ell\le\frac{t^2}{2}.
\]
Hence at most half of the pairs have squared norm greater than $t^2$.
This proves the lemma.
\end{proof}

\begin{lemma}\label{tri:upd:lem:bounded-images}
For every $0<\rho\le1$, there are an integer $L\ge3$ and $n_0$ with
the following property. Put $t=\rho/(8L)$. Let $P$ be proper,
$|P|\le n$, and $n\ge n_0$. Suppose $v:\mathbb R^2\to[0,1]$ is
finitely supported and
\[
 \sum_xv(x)\le n,\qquad \sum_{x\in P}v(x)\ge\rho n.
\]
Then, for every finitely supported real-valued function $w$, at most $L$ direct
isometries $g$ satisfy $\|gv-w\|_n\le t$.
\end{lemma}
\begin{proof}
Let $r_3(L)$ be the largest size of a three-term-progression-free
subset of $\{0,\ldots,L-1\}$. By Roth's theorem~\cite{Roth}, choose $L\ge3$ so that
$r_3(L)/L\le\rho/2$, and choose $n_0>4/\rho$.
If more than $L$ such motions existed, two distinct ones, $g,g'$,
would have rotation angles in the same one of $L$ equal half-open
intervals partitioning $[0,2\pi)$. Set $f=g^{-1}g'$.
If the two angles coincide, $f$ is a nonzero translation because
$g\ne g'$. Otherwise, $f$ is a rotation through an angle $\theta$
with $0<|\theta|<2\pi/L$. Norm invariance and telescoping give
\[
 \|fv-v\|_n=\|g'v-gv\|_n\le2t,
 \qquad \|f^jv-v\|_n\le j\|fv-v\|_n\le2jt\le\rho/4
 \quad(0\le j<L).
\]
Since $\|\mathbf1_P\|_n\le1$, the Cauchy--Schwarz inequality gives
\[
 \begin{aligned}
 \frac1n\sum_{x\in P}(f^jv)(x)
 &\ge\frac1n\sum_{x\in P}v(x)
       -\|\mathbf1_P\|_n\|f^jv-v\|_n\\
 &\ge\rho-\frac\rho4=\frac{3\rho}{4}
 \qquad(0\le j<L).
 \end{aligned}
\]
For each nonfixed $x$, the first $L$ iterates are distinct.
This is immediate for a nonzero translation. For a rotation, a
repetition would require $k\theta\in2\pi\mathbb Z$ for some
$1\le k<L$, whereas $0<|k\theta|<2\pi$.
The visit times to $P$ cannot contain $i,i+d,i+2d$ with $d>0$.
Such times would give three distinct points of $P$ satisfying
\[
 |f^ix-f^{i+d}x|=|f^{i+d}x-f^{i+2d}x|,
\]
because $f^d$ maps the first segment to the second. This contradicts
properness. Thus, each nonfixed $x$ has at most $r_3(L)$ visits.

To count these visits with their weights, substitute $y=f^jx$ to obtain
\[
 \sum_{y\in P}(f^jv)(y)
 =\sum_y\mathbf1_P(y)v(f^{-j}y)
 =\sum_xv(x)\mathbf1_P(f^jx).
\]
Summing over $j$ therefore counts each visit of the orbit starting at
$x$ with weight $v(x)$. There is at most one fixed point, and its
weight is at most one. Since $\sum_xv(x)\le n$, we obtain
\[
 \begin{aligned}
 \frac1L\sum_{j=0}^{L-1}\frac1n\sum_{y\in P}(f^jv)(y)
 &=\frac1n\sum_xv(x)
       \frac{|\{0\le j<L:f^jx\in P\}|}{L}\\
 &\le\frac{r_3(L)}L+\frac1n
 \le\frac\rho2+\frac1n.
 \end{aligned}
\]
For $n\ge n_0>4/\rho$, this contradicts the lower bound $3\rho/4$.
This proves the lemma.
\end{proof}

\begin{proof}[Proof of Theorem~\ref{tri:upd:thm:rich}]
We prove~\eqref{tri:upd:eq:rich-power-saving}. If $\varepsilon>1$,
there are no such isometries, since $|P|\le n$. Thus, assume
$0<\varepsilon\le1$. Let $P$ be proper and let $Q$ be any planar
set, both of size at most $n$. Put $\rho=\varepsilon^2/2$, and
choose $L,n_0$ and $t$ from Lemma~\ref{tri:upd:lem:bounded-images}. Fix
$\ell=\lceil4/t^2\rceil$; all these parameters depend only on
$\varepsilon$. Let $\mathcal S$ be the family of direct isometries with at least
$\varepsilon n$ matches between $P$ and $Q$. Take $n\ge n_0$. Every member then has at least two matches, and a direct isometry
is determined by the images of two distinct points.
Thus, $\mathcal S$ is finite. Write $m=|\mathcal S|$.
The assertion is immediate if $m=0$, so assume $m>0$.

Set
\[
 \mathcal D=\{f\in\mathrm{SE}(2):|P\cap fP|\ge\rho n\},
 \qquad N_g=\mathcal S\cap g\mathcal D\quad(g\in\mathcal S).
\]
The elementary estimate~\eqref{tri:upd:eq:self-rich} gives
$|\mathcal D|\le 2\rho^{-2}n$. Moreover,
$g'\in N_g$ if and only if $g'\in\mathcal S$ and
$|g'P\cap gP|\ge\rho n$.
The Cauchy--Schwarz inequality on $Q$ gives
\[
 \begin{aligned}
 \sum_{g,g'\in\mathcal S}|g'P\cap gP|
 &\ge\sum_{q\in Q}
       \left(\sum_{g\in\mathcal S}\mathbf1_{gP}(q)\right)^2\\
 &\ge\frac{(\varepsilon nm)^2}{|Q|}\ge2\rho nm^2.
 \end{aligned}
\]
The pairs with $g'\notin N_g$ contribute at most $\rho nm^2$, whereas
every remaining pair contributes at most $n$. Consequently,
\[
 \sum_{g\in\mathcal S}|N_g|\ge\rho m^2.
\]

All tuples below are ordered and allow repetitions.
Write $\mathbf g=(g_1,\ldots,g_\ell)$ and
$\mathbf f=(f_1,\ldots,f_\ell)$.
We count triples $(g,\mathbf g,\mathbf f)$, where $g\in\mathcal S$,
$\mathbf g\in\mathcal S^\ell$, $\mathbf f\in\mathcal D^\ell$, and
\[
 \left\|\frac1\ell\sum_{i=1}^{\ell}\mathbf1_{g_iP}
       -g\left(\frac1\ell\sum_{i=1}^{\ell}\mathbf1_{f_iP}\right)\right\|_n
 \le t.
\]
For each $g$ with $N_g\ne\varnothing$, apply Lemma~\ref{tri:upd:lem:averages}
to the vectors $\mathbf1_{g'P}$, $g'\in N_g$, in the normalized norm.
Their squared norms equal $|P|/n\le1$. The lemma supplies at least $|N_g|^{2\ell}/2$ pairs of tuples in $N_g^\ell$
whose averages differ by at most $t$. Writing the second tuple as
$(gf_1,\ldots,gf_\ell)$, with every $f_i\in\mathcal D$, gives
distinct triples of the required kind.
By the power-mean inequality, the number of triples is therefore at least
\[
 \frac12\sum_{g\in\mathcal S}|N_g|^{2\ell}
 \ge\frac m2\left(\frac1m\sum_{g\in\mathcal S}|N_g|\right)^{2\ell}
 \ge\frac12\rho^{2\ell}m^{2\ell+1}.
\]
Conversely, fix $(\mathbf g,\mathbf f)$ and put
\[
 w=\frac1\ell\sum_{i=1}^{\ell}\mathbf1_{g_iP},\qquad
 v=\frac1\ell\sum_{i=1}^{\ell}\mathbf1_{f_iP}.
\]
Since every $f_i\in\mathcal D$, we have
\[
 0\le v\le1,\qquad \sum_xv(x)=|P|\le n,\qquad
 \sum_{x\in P}v(x)=\frac1\ell\sum_{i=1}^{\ell}|P\cap f_iP|\ge\rho n.
\]
Lemma~\ref{tri:upd:lem:bounded-images} therefore shows that the fixed tuple pair occurs for
at most $L$ choices of $g$. There are $m^\ell|\mathcal D|^\ell$ tuple pairs,
so
\[
 \frac12\rho^{2\ell}m^{2\ell+1}
 \le Lm^\ell|\mathcal D|^\ell
 \le Lm^\ell(2\rho^{-2}n)^\ell.
\]
Rearranging gives
\[
 m^{\ell+1}\le
 2L\left(\frac{2}{\rho^4}\right)^\ell n^\ell.
\]
This proves~\eqref{tri:upd:eq:rich-power-saving} with
\[
 c_\varepsilon=\frac1{\ell+1},\qquad
 K_\varepsilon=\left(2L\left(\frac2{\rho^4}\right)^\ell
                  \right)^{1/(\ell+1)}.
\]
All parameters depend only on $\varepsilon$, so the estimate is
uniform over $P,Q$.

For orientation-reversing isometries, fix a reflection $s$.
The map $g\mapsto sg$ is injective, $sg$ preserves orientation,
and
\[
 |sQ\cap sgP|=|Q\cap gP|.
\]
Applying the direct estimate to $P,sQ$ gives the same bound.
Since $c_\varepsilon>0$, the bound is $o(n)$ for each fixed
$\varepsilon$. This proves the theorem.
\end{proof}

The proof of~\eqref{tri:upd:eq:rich-power-saving} uses properness only
of $P$. The exponent $c_\varepsilon$ depends on $\varepsilon$;
the argument does not provide an exponent independent of the
richness proportion.

\subsection{Distance energy with targets on conics}
\label{sec:curve-targets}

We estimate the mixed distance energy when the target points lie on a
conic. This estimate will control the energy lost when removing
concentrations on conics in the proof of Proposition~\ref{prop:annulus}.
A \emph{real plane conic} is the zero set in $\R^2$ of a real
polynomial of degree two. It is \emph{irreducible over $\R$} if its
defining polynomial cannot be factored into two nonconstant real
affine linear polynomials. 
Related distinct-distance estimates between a circle and
a generic point set were obtained by McDonald, McDonald,
Passant, and Sahay~\cite{McDonaldEtAl}.
Their genericity assumption is different from properness
and does not directly imply our mixed-energy estimate.

\begin{theorem}\label{thm:conic-target-energy}
Let $P\subset\R^2$ be a finite proper set, and let $T\subset P$
lie on an irreducible real plane conic with infinitely many real points.
Then
\[
 E(P,T)=o(|P|^3)
\]
uniformly over $P,T$, and the conic as $|P|\to\infty$.
\end{theorem}

This theorem will be used in the proof of Proposition~\ref{prop:annulus} to bound $E(P, T)$ where $T\subset P$ is a deleted subset lying on a conic.

In the proof of Theorem \ref{thm:conic-target-energy}, we separate pin pairs related by a Euclidean symmetry
of the conic from the remaining pairs. When the conic is not a
circle, it has only finitely many symmetries, which bounds the contribution
of these pairs. For a circle, Theorem~\ref{tri:upd:thm:rich} gives the
required bound. We estimate the remaining contribution by counting
incidences between points and curves.

A quadratic rational parametrization turns the remaining distance
equations into plane curves of degree at most eight. We will apply
the following consequence of the Pach--Sharir incidence bound to
their irreducible components. For a finite point set $\mathcal Q$ and
a finite family of curves $\Gamma$, let the \emph{incidence count} be
\[
 \mathcal I(\mathcal Q,\Gamma)
 =|\{(z,\gamma)\in\mathcal Q\times\Gamma:z\in\gamma\}|.
\]

\begin{lemma}\label{lem:bounded-degree-incidences}
Let $\mathcal Q\subset\R^2$ be finite, and let $\Gamma$ be a finite
family of distinct irreducible real plane algebraic curves of degree
at most eight, each with infinitely many real points. Then
\[
 \mathcal I(\mathcal Q,\Gamma)
 \ll |\mathcal Q|^{65/129}|\Gamma|^{128/129}
       +|\mathcal Q|+|\Gamma|,
\]
with an absolute implicit constant.
\end{lemma}

\begin{proof}
By B\'ezout's theorem, any two distinct curves in $\Gamma$
intersect in at most $64$ points. It follows that any $65$
distinct points can lie on at most one curve in $\Gamma$.

Each curve can be decomposed into a bounded number of simple
arcs after removing a bounded number of cut points, where
the bounds depend only on its degree. The removed points
contribute at most $O(|\Gamma|)$ incidences.
The remaining arcs are disjoint when they come from the
same curve, while arcs from distinct curves intersect in
at most 64 points. Moreover, any 65 distinct points can
lie on at most one arc.

We may therefore apply the Pach--Sharir incidence
theorem~\cite{PachSharir} with $65$ degrees of freedom
and multiplicity $64$. This gives
\[
 \mathcal I(\mathcal Q,\Gamma)
 \ll |\mathcal Q|^{65/129}|\Gamma|^{128/129}
       +|\mathcal Q|+|\Gamma|.
\]
Since the degrees of the curves are bounded by eight,
all implicit constants are absolute. This completes the proof.
\end{proof}

\begin{proof}[Proof of Theorem~\ref{thm:conic-target-energy}]
Let $\mathcal C$ be the supporting conic.
If $T$ is empty, then $E(P,T)=0$. If $T$ consists of
a single point, properness gives $E(P,T)=|P|$.
Thus, we may assume that $|T|\ge2$.

We also observe that if $\mathcal C$ is a circle,
its center cannot belong to $P$. Otherwise, the center
would have the same positive distance to two distinct
points of $T$, contradicting properness.

For each $a,b\in P$, define the \emph{distance locus}
\[
 \Gamma_{ab}
 =\{(x,y)\in\mathcal C\times\mathcal C:
                   |x-a|^2=|y-b|^2\}.
\]
By the definition of mixed distance energy,
\[
 E(P,T)=\sum_{a,b\in P}
             |\Gamma_{ab}\cap(T\times T)|.
\]

We divide the pairs $(a,b)$ into two classes.
Call a pair \emph{exceptional} if there is a Euclidean
isometry preserving $\mathcal C$ that maps $a$ to $b$.
We first show that the exceptional pairs contribute
$o(|P|^3)$, and then estimate the remaining contribution
using incidences between points and algebraic curves.

\medskip
\noindent
\emph{Exceptional pairs.}
Suppose first that $\mathcal C$ is not a circle.
By Pach and de Zeeuw~\cite[Lemma~2.5]{PdZ},
the conic has at most eight Euclidean symmetries.
Since each symmetry maps any given $a\in P$
to at most one point $b$, there are at most
$8|P|$ exceptional pairs.

For a fixed pair $(a,b)$ and a point $x\in T$,
properness at $b$ allows at most one $y\in T$
satisfying
\[
 |x-a|^2=|y-b|^2.
\]
This remains true when the common distance is zero.
Consequently, each exceptional pair contributes
at most $|T|$ to the mixed energy. The total
contribution is therefore bounded by
\[
 8|P||T|\le8|P|^2=o(|P|^3).
\]

Now suppose that $\mathcal C$ is a circle.
By applying a similarity, we may assume that
$\mathcal C=S^1$. This normalization preserves
properness, mixed distance energy, and the
exceptional-pair condition.

If $(a,b)$ is exceptional, then $|a|=|b|$,
since every isometry preserving $S^1$ also
preserves distances from the origin.
Consider a quadruple $(a,b,x,y)$ counted by $E(P,T)$
whose pin pair $(a,b)$ is exceptional, so that
\[
 |a|=|b|,\qquad |x|=|y|=1,\qquad
 |a-x|^2=|b-y|^2.
\]
Expanding the last equality gives
\[
 a\cdot x=b\cdot y.
\]

Thus, the ordered vector pairs $(a,x)$ and $(b,y)$
have the same Gram matrix. Hence, there exists
an orthogonal map $U\in\mathrm O(2)$ such that
\[
 Ua=b,\qquad Ux=y.
\]
This remains true when the vectors are linearly
dependent, since the isometry between their spans
can be extended to an orthogonal map of the plane.

For a a finite set $A\subset\R^2$, write
\[
 \kappa_A(U)=|A\cap UA|.
\]
For each fixed $U$, the product
$\kappa_P(U)\kappa_T(U)$ counts the quadruples
$(a,b,x,y)\in P^2\times T^2$ satisfying
$Ua=b$ and $Ux=y$. Therefore, the exceptional
contribution is at most
\[
 \sum_{U\in\mathrm O(2)}
       \kappa_P(U)\kappa_T(U).
\]

For every ordered pair $x,y\in T$, there are
exactly two orthogonal maps sending $x$ to $y$,
one preserving orientation and one reversing it.
Consequently, only finitely many summands are
nonzero, and
\[
 \sum_{U\in\mathrm O(2)}\kappa_T(U)=2|T|^2.
\]

Fix $\varepsilon>0$ and divide the sum according
to whether $\kappa_P(U)<\varepsilon|P|$ or
$\kappa_P(U)\ge\varepsilon|P|$.

For the first class, we have
\[
 \begin{aligned}
 \sum_{\substack{U\in\mathrm O(2)\\
                 \kappa_P(U)<\varepsilon|P|}}
       \kappa_P(U)\kappa_T(U)
 &\le \varepsilon|P|
       \sum_{U\in\mathrm O(2)}\kappa_T(U)\\
 &=2\varepsilon|P||T|^2
 \le2\varepsilon|P|^3.
 \end{aligned}
\]

For the second class, Theorem~\ref{tri:upd:thm:rich},
applied with both point sets equal to $P$ and
including orientation-reversing isometries,
shows that there are only $o_\varepsilon(|P|)$
orthogonal maps satisfying
$\kappa_P(U)\ge\varepsilon|P|$.
Since
\[
 \kappa_P(U)\kappa_T(U)\le|P|^2,
\]
these maps contribute at most
$o_\varepsilon(|P|^3)$.

Combining the two estimates, we obtain
\[
 \sum_{U\in\mathrm O(2)}
       \kappa_P(U)\kappa_T(U)
 \le2\varepsilon|P|^3
       +o_\varepsilon(|P|^3).
\]
Dividing by $|P|^3$, letting $|P|\to\infty$,
and then letting $\varepsilon\downarrow0$,
we conclude that the exceptional contribution
is $o(|P|^3)$, uniformly over the circle and
the point sets.

\medskip
\noindent
\emph{Intersections of the remaining distance loci.}
We now consider pairs $(a,b)$ that are not exceptional.
The key observation is that no three distance loci
corresponding to distinct nonexceptional pairs can
have infinitely many common real points.
This will allow us to control the multiplicity
of the algebraic curves arising in the incidence
argument.

We first note that every coordinate fiber of
$\Gamma_{ab}$ contains at most four points.
Indeed, if $x\in\mathcal C$ is fixed, then
\[
 |y-b|^2=|x-a|^2
\]
confines $y$ to the intersection of $\mathcal C$
with a circle centered at $b$, or to the singleton
$\{b\}$ when the distance is zero.
Since $\mathcal C$ cannot be a circle centered
at a point of $P$, B\'ezout's theorem gives at
most four possibilities for $y$.
The same argument applies when $y$ is fixed.

Suppose, for a contradiction, that three loci
\[
 \Gamma_{a_1b_1},\qquad
 \Gamma_{a_2b_2},\qquad
 \Gamma_{a_3b_3},
\]
corresponding to distinct nonexceptional pairs,
have infinitely many common real points.
By the fiber bound, their common set has
infinite projections onto both copies of
$\mathcal C$.

We first claim that the points $a_1,a_2,a_3$
are pairwise distinct, and so are $b_1,b_2,b_3$.
Indeed, suppose that $a_i=a_j$ for some $i\ne j$.
Since the pairs are distinct, we must have
$b_i\ne b_j$. Subtracting the corresponding
distance equations gives
\[
 2(b_i-b_j)\cdot y=|b_i|^2-|b_j|^2.
\]
Thus, $y$ lies on a fixed line, which intersects
the irreducible conic $\mathcal C$ in at most
two points. This contradicts the fact that
the common set has an infinite projection
onto the second copy of $\mathcal C$.
The same argument, with $x$ and $y$ interchanged,
shows that the $b_i$ are pairwise distinct.

Next, suppose that
\[
 |a_i-a_j|=|b_i-b_j|
\]
for some $i\ne j$.
The argument of Pach and de Zeeuw
in~\cite[Lemma~4.1]{PdZ} then gives a
Euclidean symmetry of $\mathcal C$ mapping
$a_i$ to $b_i$. This contradicts the
assumption that $(a_i,b_i)$ is nonexceptional.
Therefore,
\[
 |a_i-a_j|\ne|b_i-b_j|
 \qquad(i\ne j).
\]

   We first reduce the common
distance equations to an invertible affine
map, without assuming that either triple of
pins is noncollinear.

Subtracting the equation
\[
 |x-a_1|^2=|y-b_1|^2
\]
from the corresponding equations for
$i=2,3$, we obtain
\[
 2(a_i-a_1)\cdot x-2(b_i-b_1)\cdot y
 =|a_i|^2-|a_1|^2-|b_i|^2+|b_1|^2,
 \qquad i=2,3.
\]

The linear parts of these two equations
in $(x,y)\in\R^4$ are independent.
Otherwise, since the $a_i$ and $b_i$ are
pairwise distinct, there would exist
$\lambda\notin\{0,1\}$ such that
\[
 a_3-a_1=\lambda(a_2-a_1),\qquad
 b_3-b_1=\lambda(b_2-b_1).
\]
For the two affine equations to have a common
solution, their constant terms must satisfy
the same linear relation. Substituting the
displayed identities then gives
\[
 \lambda(\lambda-1)
 \bigl(|a_2-a_1|^2-|b_2-b_1|^2\bigr)=0,
\]
contrary to the pairwise distance inequality
established above.

It follows that the common solution set of
the two affine equations is a two-dimensional
affine plane in $\R^4$.

Both coordinate projections of this plane
onto the $x$- and $y$-planes must have rank two.
Indeed, if either projection had rank at most
one, its image would be contained in a point
or a line. Such a set intersects $\mathcal C$
in at most two points, contradicting the
infinite projection of the common distance loci.

Consequently, the affine plane is the graph
of an invertible affine map
\[
 y=L(x).
\]
The common distance loci contain infinitely
many points on this graph with
$x,y\in\mathcal C$. Hence, $L(\mathcal C)$
and $\mathcal C$ have infinitely many
common points. By B\'ezout's theorem,
\[
 L(\mathcal C)=\mathcal C.
\]

Moreover, the three original distance
equations hold at infinitely many points
of this graph.
We treat the parabola case directly. By applying the same
Euclidean similarity to the conic and all pins, we may assume that
\[
 \mathcal C=\{(t,t^2):t\in\R\}.
\]
An invertible affine map preserving $\mathcal C$ has the form
\[
 L(x,y)=
 \bigl(\lambda x+\mu,\,
       2\lambda\mu x+\lambda^2y+\mu^2\bigr),
 \qquad \lambda\ne0.
\]
To see this, write
$L(x,y)=(Ax+By+e,Cx+Dy+f)$ and substitute $(t,t^2)$ into
$L(\mathcal C)=\mathcal C$. Comparing coefficients gives
$B=0$, $D=A^2$, $C=2Ae$, and $f=e^2$; invertibility gives $A\ne0$.

For each $i$, the equality
\[
 |(t,t^2)-a_i|^2
 =|(\lambda t+\mu,(\lambda t+\mu)^2)-b_i|^2
\]
holds for infinitely many $t$, and hence identically.
Comparing the coefficients of $t^4$ and $t^3$ yields
\[
 \lambda^4=1,\qquad 4\lambda^3\mu=0.
\]
Thus $\lambda=\pm1$ and $\mu=0$, so $L$ is a Euclidean
symmetry of $\mathcal C$. The distance identity now gives
\[
 |x-a_i|=|x-L^{-1}(b_i)|
\]
for infinitely many $x\in\mathcal C$.
If $a_i\ne L^{-1}(b_i)$, these points would lie on a perpendicular
bisector, contradicting the irreducibility of the conic.
Therefore $b_i=L(a_i)$, contrary to the nonexceptionality of the pairs.

For an ellipse or a hyperbola, translate the center of
$\mathcal C$ to the origin and write
\[
 \mathcal C=\{x\in\R^2:x^TAx=1\},
\]
where $A$ is an invertible symmetric matrix.
Since an affine self-map of $\mathcal C$ preserves its
center, we have $L(x)=Mx$ for some invertible matrix $M$,
with
\[
 M^TAM=A.
\]
For each $i$, the identity
\[
 |x-a_i|^2=|Mx-b_i|^2
\]
holds for infinitely many $x\in\mathcal C$.
The difference of the two sides is a quadratic
polynomial vanishing on the irreducible conic
$\mathcal C$, so it is a scalar multiple of
$x^TAx-1$. Comparing the quadratic and linear
terms gives
\[
 I-M^TM=cA,\qquad a_i=M^Tb_i
\]
for some $c\in\R$.
On the other hand, $M^TAM=A$ implies
$MA^{-1}M^T=A^{-1}$. Taking traces, we obtain
\[
\begin{aligned}
 2c
 &=\operatorname{tr}\bigl((I-M^TM)A^{-1}\bigr)\\
 &=\operatorname{tr}(A^{-1})
   -\operatorname{tr}(MA^{-1}M^T)
 =0.
\end{aligned}
\]
Thus $M^TM=I$, so $M$ is orthogonal, and
$b_i=Ma_i=L(a_i)$. This makes $(a_i,b_i)$
exceptional, again a contradiction.

This contradiction proves that no three
distinct nonexceptional distance loci
can have infinitely many common real points.

\needspace{5\baselineskip}
\emph{Parametrization and the incidence bound.}
An irreducible real conic with infinitely many real points has a
nonsingular projective closure and a quadratic rational
parametrization over $\R$; see \cite[Problem~6.37]{Fulton}.
Choose an affine parameter chart whose omitted point does not
correspond to a point of $T$. In this chart write
\[
 \phi(u)=\left(\frac{X(u)}{Z(u)},\frac{Y(u)}{Z(u)}\right),
 \qquad \deg X,\deg Y,\deg Z\le2.
\]
The parametrization is one-to-one on its domain $Z(u)\ne0$ and
covers $T$. Thus, $\mathcal T=\phi^{-1}(T)\subset\R$ has
$|\mathcal T|=|T|$, and $Z$ is nonzero on $\mathcal T$.

For $a=(a_x,a_y)\in P$, define
\[
 Q_a(u)=(X(u)-a_xZ(u))^2+(Y(u)-a_yZ(u))^2.
\]
Then
\[
 |\phi(u)-a|^2=\frac{Q_a(u)}{Z(u)^2}
 \qquad\text{when }Z(u)\ne0.
\]
For each nonexceptional pair $(a,b)$, clearing denominators in
$|\phi(u)-a|^2=|\phi(v)-b|^2$ gives the polynomial
\[
 F_{ab}(u,v)=Q_a(u)Z(v)^2-Q_b(v)Z(u)^2.
\]
Since $\deg Q_a,\deg Q_b\le4$ and $\deg Z^2\le4$, the total
degree of $F_{ab}$ is at most eight. Its zeros on $\mathcal T^2$
correspond exactly to
$\Gamma_{ab}\cap(T\times T)$.
It is not the zero polynomial: otherwise varying $u$ and $v$
independently would make the distance from $a$ constant on
$\mathcal C$, forcing $\mathcal C$ to be a circle centered at $a$,
which has already been excluded.

Factor each $F_{ab}$ over $\R$, ignoring repeated factors.
Discard factors dividing $Z(u)Z(v)$, since they vanish at no point
of $\mathcal T^2$.
Also discard factors with only finitely many real zeros. Each has
only a bounded number of such zeros, since a real algebraic curve
of bounded degree has a bounded number of connected components;
see \cite[Theorem~2.2]{PdZ}.
Their total contribution is $O(|P|^2)$.

Every remaining irreducible factor has infinitely many real zeros
off $Z(u)Z(v)=0$, by B\'ezout. If it occurred for three distinct
nonexceptional pin pairs, injectivity of $\phi\times\phi$ would give
infinitely many common real points of the three corresponding
distance loci, contradicting the claim above. Thus, each retained
curve occurs for at most two pin pairs.

Let $\Gamma$ be the family of distinct retained curves. There are
at most $8|P|^2$ of them, and each has degree at most eight.
Accounting for their multiplicity and the discarded finite sets,
Lemma~\ref{lem:bounded-degree-incidences} bounds the
nonexceptional contribution by
\[
 \begin{aligned}
 2\mathcal I(\mathcal T^2,\Gamma)+O(|P|^2)
 &\ll (|T|^2)^{65/129}(|P|^2)^{128/129}
       +|T|^2+|P|^2\\
 &\ll |P|^{3-1/129},
 \end{aligned}
\]
where the last inequality uses $|T|\le|P|$.
All constants are absolute, since the degrees and component
multiplicities are bounded independently of the conic.
Adding the exceptional contribution proves
$E(P,T)=o(|P|^3)$ with the required uniformity. This completes the proof.
\end{proof}

\begin{remark}
The estimate in Theorem~\ref{thm:conic-target-energy}
cannot in general be improved to $O(|P|^{3-\delta})$
for any fixed $\delta>0$, even when the supporting
conic is the unit circle. Indeed, let $A\subset\{1,\ldots,M\}$ be a Behrend set
containing no nontrivial three-term arithmetic progression,
with
\[
 n=|A|\ge M\exp(-C\sqrt{\log M}).
\]
Set
\[
 P=\left\{e^{2\pi ia/(4M)}:a\in A\right\}
 \subset S^1.
\]
The absence of three-term arithmetic progressions
implies that $P$ is proper. Since the chord lengths
depend only on the differences of the indices,
we have $D(P)\le M-1$. By Cauchy--Schwarz,
\[
 E(P)\ge\frac{n^2(n-1)^2}{M-1}
 \gg n^3\exp(-C'\sqrt{\log n})
 =n^{3-o(1)}.
\]
Since $E(P,P)=E(P)+n^2$, the same obstruction
applies to the mixed-energy estimate.
\end{remark}

\subsection{Distance energy in shrinking annuli}\label{sec:annulus-proof}
We prove Proposition~\ref{prop:annulus} using the motion and conic
estimates from Subsections~\ref{tri:upd:sec:rich}
and~\ref{sec:curve-targets}. For $0<t\le2\pi$, define
\[
 \Omega_j(t)=\frac1{|P_j|}\sup_{|I|\le t}
 \bigl|\{x\in P'_j:\arg x\in I\}\bigr|.
\]

For every $x\in P'_j$, we have $|x|\ge1-h_j>1/2$, so $x\ne0$
and $\arg x$ is well defined. The quantity $\Omega_j(t)$ is the
largest number of points of $P'_j$ whose angles lie in an angular
interval of length at most $t$, divided by $|P_j|$.
In particular, $\Omega_j$ is nonnegative and nondecreasing in $t$.
The \emph{angular nonconcentration} hypothesis
\eqref{eq:angular-nonconcentration} therefore has the equivalent form
\par
\[
 \lim_{t\downarrow0}\limsup_{j\to\infty}\Omega_j(t)=0.
\]

Indeed, for every $\varepsilon>0$, the hypothesis supplies $t_0>0$
with $\limsup_j\Omega_j(t_0)\le\varepsilon$.
Monotonicity then gives the displayed limit.
Conversely, choose $t_0$ with $\limsup_j\Omega_j(t_0)<\varepsilon$;
then $\Omega_j(t_0)\le\varepsilon$ for all sufficiently large $j$.
We will also use the consequence
\begin{equation}\label{eq:shrinking-angular-intervals}
 \Omega_j(t_j)\longrightarrow0
 \quad\text{for every positive sequence }t_j\longrightarrow0.
\end{equation}

For each fixed $t>0$, we have $t_j\le t$ for all sufficiently large
$j$, and hence
\[
 0\le\limsup_{j\to\infty}\Omega_j(t_j)
 \le\limsup_{j\to\infty}\Omega_j(t).
\]
Letting $t\downarrow0$ proves~\eqref{eq:shrinking-angular-intervals}.
\par

All limits in this subsection are along the sequence in
Proposition~\ref{prop:annulus}. Write $P=P_j$, $P'=P'_j$, $n=|P|$,
and $h=h_j$. Angles label the point pairs in the incidence argument,
while all equal-distance conditions use their exact Euclidean distances.
It suffices to bound $E(P_j)$, since \eqref{eq:cs-distances} then
implies the bound for $D(P_j)$.

After removing concentrations on conics, we discard equal-distance
quadruples whose source angular separation is close to $0$ or $\pi$.
Each remaining quadruple is realized by two Euclidean isometries.
One nearly fixes the origin and shows that the two associated lines
have close angular labels.
The other moves the origin by a definite amount. We bound the
contribution of quadruples for which the latter isometry has
at least three matches;
a line-incidence estimate handles
those for which it has exactly two matches.

\subsubsection{Removing concentrations on conics}
Fix $\sigma>0$ independently of $j$ and $n$.
We keep $\sigma$ fixed while letting $j\to\infty$, and let
$\sigma\downarrow0$ only at the end of the proof.
We first remove concentrations on conics while losing only
$o_\sigma(n^3)$ equal-distance quadruples.
Lemma~\ref{lem:restriction} gives
\[
 E(P)\le E(P')+4|P\setminus P'|n^2=E(P')+o(n^3).
\]

\begin{figure}[H]
\centering

\begin{tikzpicture}[
  x=2.05cm,y=2.05cm,font=\small,text=blue,
  boundary/.style={blue!45,thin},
  ray/.style={blue!60,densely dotted,thin},
  part/.style={blue,line width=1.8pt},
  anglepart/.style={blue!65!black,line width=1.8pt}
]
\pgfmathsetmacro{\figlo}{acos(.55/.9)}
\pgfmathsetmacro{\fighi}{acos(.55/1.1)}
\pgfmathsetmacro{\figylo}{sqrt(.9*.9-.55*.55)}
\pgfmathsetmacro{\figyhi}{sqrt(1.1*1.1-.55*.55)}
\pgfmathsetmacro{\figtangent}{acos(1/1.1)}
\pgfmathsetmacro{\figytangent}{sqrt(1.1*1.1-1)}
\foreach \shift in {0,3.45}{
  \begin{scope}[shift={(\shift,0)}]
    \fill[blue!6,even odd rule] (0,0) circle (1.1) (0,0) circle (.9);
    \draw[boundary] (0,0) circle (.9);
    \draw[boundary] (0,0) circle (1.1);
    \draw[blue!60,dashed,thin] (0,0) circle (1);
    \fill[blue] (0,0) circle (.018);
    \node[below left,inner sep=2pt] at (0,0) {$O$};
    \draw[blue,-{Latex[length=1.4mm]}] (0,0) -- (.38,0);
    \node[below,inner sep=2pt] at (.25,0) {$u$};
  \end{scope}
}

\node at (0,1.47) {(a) Two angular intervals};
\draw[blue!65,thin] (.55,-1.22) -- (.55,1.22);
\node[above right,inner sep=2pt] at (.55,1.16) {$\ell$};
\draw[part] (.55,\figylo) -- (.55,\figyhi);
\draw[part] (.55,-\figylo) -- (.55,-\figyhi);
\draw[anglepart] (\figlo:1)
  arc[start angle=\figlo,end angle=\fighi,radius=1];
\draw[anglepart] (-\fighi:1)
  arc[start angle=-\fighi,end angle=-\figlo,radius=1];
\foreach \ang in {\figlo,\fighi,-\figlo,-\fighi}
  \draw[ray] (0,0) -- (\ang:1.18);
\node at (.19,1.10) {$I_1$};
\node at (.19,-1.10) {$I_2$};
\draw[boundary] (145:1.1) -- (-1.14,.89);
\node[left,inner sep=1pt] at (-1.14,.89) {$1+h$};
\draw[boundary] (155:.9) -- (-1.13,.41);
\node[left,inner sep=1pt] at (-1.13,.41) {$1-h$};
\begin{scope}[shift={(3.45,0)}]
  \node at (0,1.47) {(b) Tangency to the unit circle};
  \draw[blue!65,thin] (1,-1.22) -- (1,1.22);
  \draw[part] (1,-\figytangent) -- (1,\figytangent);
  \draw[anglepart] (-\figtangent:1)
    arc[start angle=-\figtangent,end angle=\figtangent,radius=1];
  \draw[ray] (0,0) -- (-\figtangent:1.19);
  \draw[blue,thin] (0,0) -- (1,\figytangent);
  \draw[blue,thin] (.39,0) -- (1,0);
  \draw[blue,thin] (.92,0) -- (.92,.08) -- (1,.08);
  \draw[blue,thin] (.34,0)
    arc[start angle=0,end angle=\figtangent,radius=.34];
  \node at (.48,.12) {$\alpha_h$};
  \node at (.57,.34) {$1+h$};
  \node[below,inner sep=2pt] at (.68,0) {$1$};
  \fill[blue] (1,\figytangent) circle (.018);
  \node[above right,inner sep=2pt] at (1,\figytangent) {$Q$};
  \node[below right,inner sep=2pt] at (1,0) {$T$};
  \node at (1.22,-.34) {$I$};
\end{scope}
\end{tikzpicture}
\caption{The thick line segments lie in the annulus;
radial projection places their angles in the highlighted arcs.
In the tangent example, $OT=1$ and $OQ=1+h$, so
$\cos\alpha_h=1/(1+h)$. The angular interval has length
$2\alpha_h\sim2\sqrt{2h}$ as $h\downarrow0$, which explains the
square-root scale in the uniform bound.}
\label{fig:annulus-angular-cover}
\end{figure}
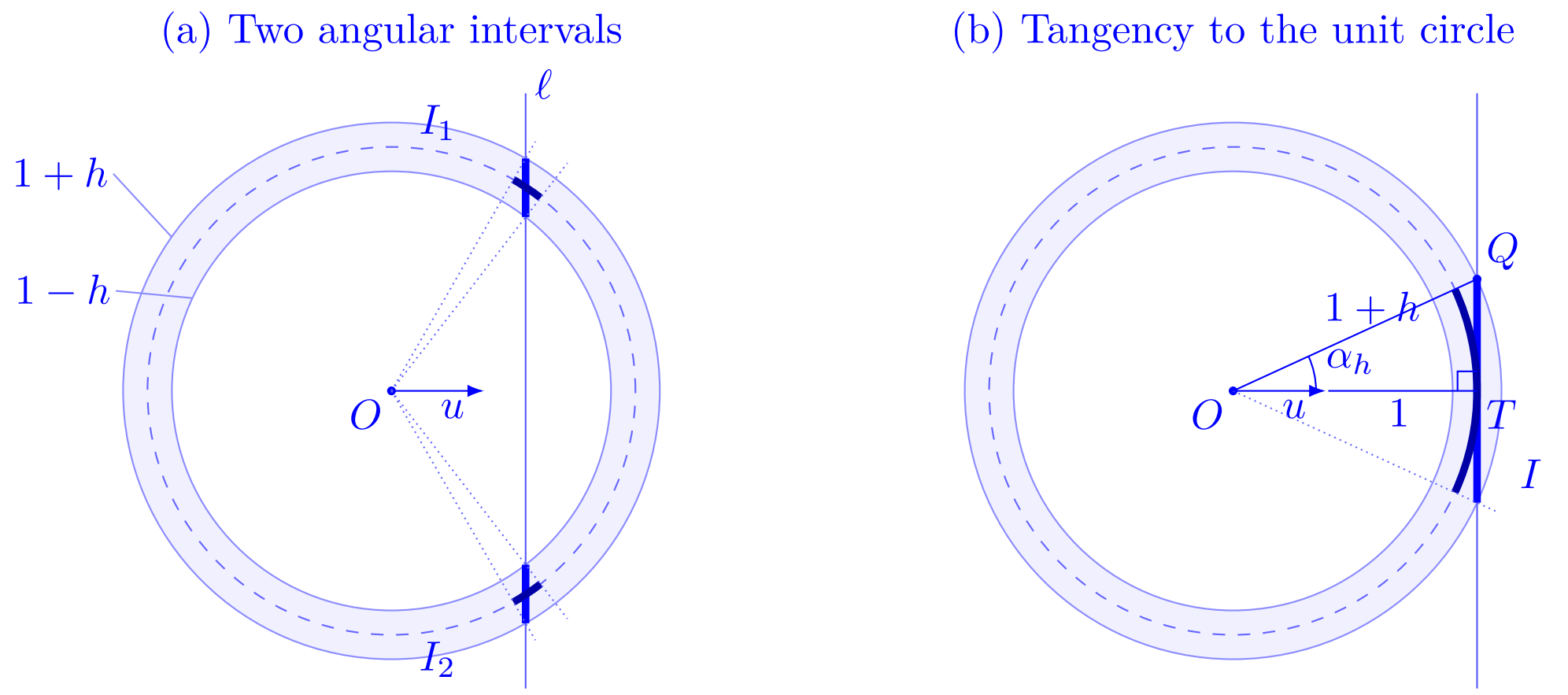

We first bound the number of points of $P'$ on any line.
For a line $\ell$, the angles of the points in its intersection
with the annulus lie in at most two angular intervals of length
$O(\sqrt h)$; see Figure~\ref{fig:annulus-angular-cover}.
To prove this, write $\ell=\{x:x\cdot u=b\}$, where $u$ is a unit
normal vector. If $x$ belongs to this intersection and
$\varphi=\arg x-\arg u$, then
\[
 |x|\cos\varphi=b.
\]
Since $1-h\le|x|\le1+h$ and $|b|\le1+h$ whenever the intersection
is nonempty, we have
\[
 |\cos\varphi-b|
 =\left|\frac{b}{|x|}-b\right|
 \le\frac{(1+h)h}{1-h}\le3h.
\]
Thus, $\cos\varphi$ belongs to an interval of length at most $6h$
in $[-1,1]$. On each of $[0,\pi]$ and $[\pi,2\pi]$, cosine is
monotone, so the inverse image of that interval is an interval
or is empty. To bound its length uniformly, including near the
endpoints, observe that for $0\le\alpha\le\beta\le\pi$,
\par
\[
 \cos \alpha-\cos \beta
 =2\sin\frac{\alpha+\beta}{2}\sin\frac{\beta-\alpha}{2}
 \ge2\sin^2\frac{\beta-\alpha}{2}
 \ge\frac{2}{\pi^2}(\beta-\alpha)^2.
\]

Here the middle inequality follows from
\[
 \frac{\beta-\alpha}{2}\le\frac{\alpha+\beta}{2}
 \le\pi-\frac{\beta-\alpha}{2},
\]
and the last uses $\sin s\ge2s/\pi$ for $0\le s\le\pi/2$.
If $[\alpha,\beta]$ is an inverse image under consideration, then
$\cos\alpha-\cos\beta\le6h$, so $\beta-\alpha\le\pi\sqrt{3h}$.
Reflection by $\varphi\mapsto2\pi-\varphi$ gives the same bound
on $[\pi,2\pi]$. This proves the angular covering with an
absolute constant $C$.

Each of the two covering intervals contains the angles of at most
$n\Omega_j(C\sqrt h)$ points of $P'$. The intervals may depend on
$\ell$, but this bound does not. Consequently,
\[
 \sup_{\ell\text{ a line}}|P'\cap\ell|
 \le2n\Omega_j(C\sqrt h)=o(n).
\]
The final estimate is along $j\to\infty$: since $h=h_j\to0$,
\eqref{eq:shrinking-angular-intervals}, with
$t_j=C\sqrt{h_j}$, gives $\Omega_j(C\sqrt{h_j})\to0$.
\par

Starting with $P'$, we repeatedly delete all remaining points on an
irreducible conic containing at least $\sigma n$ such points.
The threshold $\sigma n$ always uses the original cardinality $n=|P|$.
If $k$ steps are performed and their deleted sets are
$T_1,\ldots,T_k$, then these sets are pairwise disjoint and
\[
 k\sigma n\le\sum_{i=1}^k|T_i|\le|P'|\le n.
\]
Thus, there are at most $1/\sigma$ steps.
Let $S$ be the set remaining after these deletions.
At termination, every irreducible conic of degree two contains
fewer than $\sigma n$ points of $S$.

We next check the hypothesis on real points needed to apply
Theorem~\ref{thm:conic-target-energy} to each deleted set $T$.
An irreducible real conic containing two distinct real points
has infinitely many real points. Indeed, if $F$ defines the conic
and $z_0,z_1$ are two distinct real zeros, then
\[
 F\bigl(z_0+\lambda(z_1-z_0)\bigr)
 =c\lambda(\lambda-1),\qquad c\ne0.
\]
The restriction cannot vanish identically, since that would give
$F$ a real linear factor. Its derivative at $\lambda=0$ is
nonzero, so $\nabla F(z_0)\ne0$; the implicit function theorem
then gives a real arc of zeros near $z_0$.
Since $\sigma$ is fixed and $n=|P_j|\to\infty$, we have
$|T|\ge\sigma n\ge2$ for all sufficiently large $j$.
Consequently, every conic selected in the deletion procedure
has infinitely many real points, although each deleted set $T$
is finite.
\par

Apply Theorem~\ref{thm:conic-target-energy} to the
original proper set $P$ and its subset $T$. This gives
$E(P,T)=o(n^3)$, uniformly in the chosen conic and in the deletion step.
The number of equal-distance quadruples involving $T$ in
one specified endpoint position is at most
\[
 \sum_{r>0}\bigl|\{(x,y)\in P\times T:|x-y|^2=r\}\bigr|\nu_r(P)
 \le\sqrt{E(P,T)E(P)}
 \le n^{3/2}E(P,T)^{1/2}=o(n^3).
\]

The sum of the squares of these mixed-distance multiplicities is at
most $E(P,T)$, which also includes the zero-distance term.
The last inequality uses $E(P)\le n^3$, which
follows from properness. The $o(n^3)$ bound is uniform
over the deleted sets $T$. Since $\sigma$ is fixed, summing over
the four endpoint positions and the at most $1/\sigma$ deletion
steps gives a total loss of $o_\sigma(n^3)$.

Every degenerate real conic is contained
in the union of at most two lines or has at most one real point.
The uniform $o(n)$ line bound therefore makes its intersection with
$P'$ smaller than $\sigma n$ for all sufficiently large $j$.
For a real polynomial $F$, let
$\mathcal{Z}(F)=\{x\in\R^2:F(x)=0\}$ denote its \emph{real zero set}.
Thus, for all sufficiently large $j$, the remaining set $S$ satisfies
\begin{equation}\label{eq:conic}
 M_2(S):=\max_{\substack{0\ne F\in\R[x_1,x_2]\\ \deg F\le2}}
             |S\cap \mathcal{Z}(F)|\le\sigma n,
 \qquad E(P)\le E(S)+o_\sigma(n^3).
\end{equation}
Throughout the remaining estimates, $P$ remains the original set,
$n=|P|$, and $\Omega_j$ is computed on $P'$. In particular, neither
$n$ nor $\Omega_j$ depends on $\sigma$, even if $|S|/n$ tends to zero.
It remains to estimate the number of equal-distance quadruples in $S^4$.

\subsubsection{Two isometries for equal-distance quadruples}

Write an equal-distance quadruple in $S^4$ as
$(p_1,p_2,q_1,q_2)$.
Its source and target segments are
$(p_1,p_2)$ and
$(q_1,q_2)$, respectively; an isometry realizing it
sends $p_i$ to $q_i$ for $i=1,2$.

Fix $0<\theta<\pi/8$, independently of $j$.
For distinct $p_1,p_2\in S$, let $\gamma\in[0,\pi]$ be the smaller
angle between the vectors $p_1,p_2$ from the origin.
Call the source segment \emph{good} if
$\gamma\in[\theta,\pi-\theta]$.
Its half-length is $|p_2-p_1|/2$, while the distance of its
midpoint from the origin is $|(p_1+p_2)/2|$.
On the unit circle these quantities are $\sin(\gamma/2)$ and
$\cos(\gamma/2)$, respectively. More generally, the identities
\[
 \begin{aligned}
 \frac{|p_2-p_1|^2}{4}
 &=\frac{(|p_2|-|p_1|)^2}{4}
     +|p_1||p_2|\sin^2(\gamma/2),\\
 \left|\frac{p_1+p_2}{2}\right|^2
 &=\frac{(|p_2|-|p_1|)^2}{4}
     +|p_1||p_2|\cos^2(\gamma/2)
 \end{aligned}
\]
show that, for a good source segment in the annulus,
\[
 \min\left\{\frac{|p_2-p_1|}{2},
             \left|\frac{p_1+p_2}{2}\right|\right\}
 \ge(1-h)\sin(\theta/2)\ge\tfrac12\sin(\theta/2)>0.
\]
Thus, excluding angular separations near zero prevents the segment
from becoming arbitrarily short, and excluding separations near
$\pi$ keeps its midpoint away from the origin. Both lower bounds
are uniform in $j$ and in the chosen segment, with $\theta$ fixed.

For a fixed $p_1$, a source segment that is not good has
$\arg p_2$ in the interval of length $2\theta$ centered at
$\arg p_1$, or in the interval of length $2\theta$ centered at
$\arg p_1+\pi$, with angles taken modulo $2\pi$.
Since $S\subset P'$, each interval contains the angles of at most
$n\Omega_j(2\theta)$ possible choices of $p_2$.
There are at most $n$ choices of $p_1$, so the number of ordered
source segments that are not good is at most
\[
 n\bigl(2n\Omega_j(2\theta)\bigr)
 =2\Omega_j(2\theta)n^2.
\]
For each such source segment, there are at most $n$ choices of
the first target endpoint $q_1$. Properness of $S$ leaves at most
one choice of $q_2$ with $|q_1-q_2|=|p_1-p_2|$:

two distinct choices would give the same positive distance from
the pin $q_1$. Hence, the number of equal-distance quadruples whose
source segment is not good is at most
\[
 \bigl(2\Omega_j(2\theta)n^2\bigr)n
 =2\Omega_j(2\theta)n^3.
\]
We discard these quadruples and leave the point
set $S$ unchanged. At the end of the proof, we divide this error
by $n^3$, take the limsup as $j\to\infty$ with $\theta$ fixed,
and then let $\theta\downarrow0$; angular nonconcentration makes
the resulting bound tend to zero. We now consider only
equal-distance quadruples whose source segment is good.

For a Euclidean isometry $g$, its \emph{orthogonal part} $O_g$
is the orthogonal linear map in the decomposition
\[
 g(x)=O_gx+g(0),\qquad |O_gx|=|x|\quad(x\in\R^2).
\]
It is a rotation about the origin, including the identity, when
$\det O_g=1$, and a reflection in a line through the origin
when $\det O_g=-1$. The vector $g(0)$ is the \emph{translation part}
in this decomposition, and $|g(0)|$ measures the displacement
of the origin.
\par

\begin{lemma}\label{lem:two-isometries}
Suppose that $p_1,p_2,q_1,q_2$ lie in the annulus
$\lvert |x|-1\rvert\le h$, that
$|p_1-p_2|=|q_1-q_2|>0$, and that the source segment is good.
For all sufficiently small $h$, depending only on $\theta$, the two
Euclidean isometries taking $p_i$ to $q_i$ can be denoted by $g,f$ so
that
\begin{equation}\label{eq:nearfar}
 \begin{gathered}
 0\le |g(0)|\le C_\theta h,\qquad
 0<c_\theta\le |f(0)|\le2(1+h),\\
 \det O_g=-\det O_f.
 \end{gathered}
\end{equation}
The constants depend only on $\theta$.
\end{lemma}

\begin{proof}
We compare the images of the segment midpoints under the two possible
orthogonal maps. The target midpoint minus each image is exactly
the translation vector of the corresponding isometry. Set
\[
 a=(p_1+p_2)/2,\quad v=(p_2-p_1)/2,\quad
 b=(q_1+q_2)/2,\quad w=(q_2-q_1)/2.
\]
There are exactly two orthogonal maps $O$ taking $v$ to $w$, one of
each determinant. The corresponding isometries are
$x\mapsto Ox+b-Oa$. The identities
\[
 a\cdot v=\frac{|p_2|^2-|p_1|^2}{4},\qquad
 |a|^2+|v|^2=\frac{|p_1|^2+|p_2|^2}{2},
\]
and the corresponding identities for $b,w$, give
\[
 |v|=|w|,\qquad a\cdot v,b\cdot w=O(h),\qquad
 |a|^2,|b|^2=1-|v|^2+O(h).
\]
Since the source segment is good, both $|v|^2$ and $1-|v|^2$ are bounded
below in terms of $\theta$ when $h$ is sufficiently small. Indeed,
if the smaller source angle is $\gamma\in[\theta,\pi-\theta]$, then
\[
 |v|^2
 =\frac{|p_1|^2+|p_2|^2-2|p_1||p_2|\cos\gamma}{4}
 =\sin^2(\gamma/2)+O(h).
\]
Both lower bounds can therefore be taken to be
$\tfrac12\sin^2(\theta/2)$ for sufficiently small $h$.

Let $a_\perp,b_\perp$ be the components of $a,b$ perpendicular
to $v,w$, respectively. Then $|a-a_\perp|,|b-b_\perp|=O_\theta(h)$ and
\[
 \begin{aligned}
 |a_\perp|^2&=|a|^2-\frac{(a\cdot v)^2}{|v|^2}
 =1-|v|^2+O_\theta(h),\\
 |b_\perp|^2&=|b|^2-\frac{(b\cdot w)^2}{|w|^2}
 =1-|v|^2+O_\theta(h).
 \end{aligned}
\]
In particular, $|a_\perp|$ and $|b_\perp|$ are bounded below in
terms of $\theta$, and
\[
 \bigl||a_\perp|-|b_\perp|\bigr|
 =\frac{\bigl||a_\perp|^2-|b_\perp|^2\bigr|}
 {|a_\perp|+|b_\perp|}=O_\theta(h).
\]

Choose $O_g$ so that $O_ga_\perp$ and $b_\perp$ have the same
direction, and let $O_f$ be the other orthogonal map taking $v$ to $w$.
For the corresponding isometries,
\[
 \begin{aligned}
 |g(0)|&=|b-O_ga|=O_\theta(h),\\
 |f(0)|&=|b-O_fa|\ge |a_\perp|+|b_\perp|\ge c_\theta>0.
 \end{aligned}
\]
Their orthogonal parts have opposite determinants.
For the upper bound on $|f(0)|$,
use $f(p_1)=q_1$ and the orthogonality of $O_f$ to obtain
\[
 |f(0)|=|q_1-O_fp_1|
 \le |q_1|+|p_1|\le2(1+h).
\]
Together with the preceding bounds, this proves the lemma.
\end{proof}

The first isometry may fix the origin exactly:
if $q_i=p_i$ for $i=1,2$, the choice of $g$ is the identity.
Thus, no positive lower bound for $|g(0)|$ holds in general.

We next count pairs consisting of an equal-distance quadruple and a
Euclidean isometry realizing it. For a Euclidean isometry $f$, let
\[
 \kappa_S(f)=|S\cap fS|=|\{x\in S:f(x)\in S\}|
\]
be its number of \emph{matches}.
For a fixed $f$, choosing two distinct matched source points in
order gives exactly $\kappa_S(f)(\kappa_S(f)-1)$ quadruples
$(p_1,p_2,f(p_1),f(p_2))$ counted by $E(S)$.
Conversely, every equal-distance quadruple in $S^4$ is realized by
exactly two Euclidean isometries, one of each orientation.
Counting the pairs in these two ways gives
\[
 \sum_{f:\,\kappa_S(f)\ge2}\kappa_S(f)(\kappa_S(f)-1)=2E(S).
\]

For each fixed finite set $S$, this sum has finitely many terms:
every isometry in it realizes an equal-distance quadruple from $S^4$.

For each remaining equal-distance quadruple with a good source segment,
select the realizing isometry $f$ given by
Lemma~\ref{lem:two-isometries} with
$|f(0)|\ge c_\theta$, and denote
the other realizing isometry by $g$. The two specified matches give
$\kappa_S(f)\ge2$.
The following lemma, applied with $d=c_\theta$, bounds the
contribution of quadruples whose selected $f$ satisfies
$\kappa_S(f)\ge3$ by $o_\theta(n^3)$.
This leaves those with
$\kappa_S(f)=2$ for the line-incidence argument.
\par

\begin{lemma}\label{lem:far}
For every fixed $d>0$,
\[
 \sum_{\substack{f:\ |f(0)|\ge d\\ \kappa_S(f)\ge3}}
          \kappa_S(f)(\kappa_S(f)-1)=o_d(n^3).
\]
\end{lemma}

\begin{proof}

We show that, for each isometry in the sum, a fixed positive fraction of
its ordered pairs of distinct matched source points have small angular
separation. Angular nonconcentration will then bound their total
contribution.
Write $f(x)=O_fx+v_0$ and $z=-O_f^Tv_0$. Since
$|x-z|=|O_fx+v_0|=|f(x)|$, a matched point satisfies

\[
 |x|=1+O(h),\qquad |x-z|=1+O(h).
\]
Subtracting their squares gives $2x\cdot z-|z|^2=O(h)$.
Dividing by $2|z|$, which is legitimate
because $|z|=|v_0|\ge d>0$, yields
\[
 |z|\ge d,\qquad
 \left|x\cdot\frac{z}{|z|}-\frac{|z|}{2}\right|=O_d(h).
\]
Whenever a match exists,
$|z|\le |x|+|x-z|\le2(1+h)$. Since $|x|=1+O(h)$,
the preceding estimate therefore gives
\[
 \cos(\arg x-\arg z)=\frac{|z|}{2}+O_d(h).
\]
The cosine estimate used for lines now confines all matched source
angles to at most two intervals of length $O_d(\sqrt h)$.
Assign the points disjointly to these intervals,
with populations $u,v$ and $u+v=\kappa_S(f)$. For
$\kappa_S(f)\ge3$,
\[
 u(u-1)+v(v-1)\ge\frac{\kappa_S(f)(\kappa_S(f)-1)}3.
\]

For $\kappa_S(f)=3$, the minimum on the left is $2$, so equality is possible.
For $\kappa_S(f)\ge4$, one has
\[
 u(u-1)+v(v-1)
 =u^2+v^2-\kappa_S(f)
 \ge\frac{\kappa_S(f)(\kappa_S(f)-2)}2
 \ge\frac{\kappa_S(f)(\kappa_S(f)-1)}3.
\]

Thus, at least one third of the ordered pairs of distinct matched
source points have angular separation $O_d(\sqrt h)$.

For a suitable constant $C_d$, there are at most
$n^2\Omega_j(C_d\sqrt h)$ ordered source pairs with this
separation. Each such pair extends to at most $n$
equal-distance quadruples, and each quadruple is realized by
exactly two Euclidean isometries.
Accounting for the factor three in the preceding pair bound and
the two isometries realizing each quadruple, the sum in the statement
is at most
\[
 6n^3\Omega_j(C_d\sqrt h)=o_d(n^3).
\]
This proves the lemma.
\end{proof}

For fixed $\theta$ and all sufficiently large $j$,
Lemma~\ref{lem:two-isometries} applies to every remaining
equal-distance quadruple. Apply Lemma~\ref{lem:far} with
$d=c_\theta$ from \eqref{eq:nearfar}.
The number of quadruples assigned to any fixed $f$ is at most
$\kappa_S(f)(\kappa_S(f)-1)$. Consequently, the quadruples
whose selected $f$ has at least three matches contribute
$o_\theta(n^3)$ and may be discarded.
Every selected $f$ for a remaining quadruple still has its two
specified matches, so it remains to consider $\kappa_S(f)=2$.
No estimate for $\kappa_S(g)$ is needed in this reduction.

\subsubsection{Line representations and concentration bounds}
We represent the remaining equal-distance quadruples by intersection points
of lines. For each fixed $U$ below, the condition $\kappa_S(f)=2$
ensures that each selected point represents at most two such quadruples.
Let $\mathrm{Id}$ denote the identity map, and let $R$ be a fixed reflection
in a line through the origin. For each remaining equal-distance quadruple,
choose $U\in\{\mathrm{Id},-\mathrm{Id},R,-R\}$ so that the composition
$Uf=U\circ f$ is an orientation-preserving isometry with rotation angle
$\alpha\in[\pi/2,3\pi/2]$. Indeed, the two choices of $U$ with
determinant $\det O_f$ are negatives of one another, so the rotation
angles of the resulting orthogonal maps $UO_f$ differ by $\pi$
modulo $2\pi$. The closed interval $[\pi/2,3\pi/2]$ contains
at least one of these angles, so such a choice always exists,
including when $f$ is a translation or a reflection.
In particular, $Uf$ is a nonidentity rotation.

To count the selected rotations as intersections of lines, we use the
rotation coordinates of Guth and Katz~\cite[Proposition~2.7]{GuthKatz}.
Let $J$ be counterclockwise rotation through $\pi/2$.
For fixed $U$, define the line family
\[
 \mathcal{L}_U=\{\ell_{p,Uq}:p,q\in S\},\qquad
 \ell_{p,Uq}:\quad
 c(\tau)=(p+Uq)/2+(\tau/2)J(Uq-p),\quad \tau\in\R.
\]

In the coordinates $(c,\tau)\in\R^2\times\R$, this line represents
rotations taking $p$ to $Uq$, with center $c$ and
$\tau=\cot(\alpha/2)$. The chosen angle range gives $|\tau|\le1$, so every
selected motion is represented by a finite point. Distinct endpoint
pairs determine distinct lines, including when $p=Uq$.
The point representing $Uf$ lies on exactly $\kappa_S(f)$ lines of
$\mathcal{L}_U$, since $Uf(x)\in US$ if and only if $f(x)\in S$.
Consequently, when $\kappa_S(f)=2$, the two lines through that point
determine the two matched endpoint pairs. Their two possible orderings
give at most two equal-distance quadruples.

A \emph{regulus} is a nonplanar doubly ruled real quadric surface.
The following lemma controls the concentrations of lines that enter
the incidence bound.

\begin{lemma}\label{lem:plane-regulus-bounds}
For each fixed $U\in\{\mathrm{Id},-\mathrm{Id},R,-R\}$, one can delete $o_\sigma(n^2)$
lines from $\mathcal{L}_U$ so that every plane and every regulus contains
at most $4\sigma n$ remaining lines. These deletions remove at most
$o_\sigma(n^3)$ ordered pairs of intersecting lines.
\end{lemma}

\begin{proof}
We first delete the lines in planes containing at least $\sigma n$
lines, and bound the resulting loss. The conic concentration bound
\eqref{eq:conic} will then control every regulus without further
deletions. A nonhorizontal plane has an equation
$a\cdot c+b\tau=e$, where $|a|=1$. It contains $\ell_{p,y}$ exactly when
\[
 a\cdot(p+y)=2e,\qquad a\cdot J(y-p)=-2b.
\]
In orthonormal coordinates with axes $a,Ja$, these equations become
$y_1=2e-p_1$ and $y_2=p_2+2b$. Thus, the plane corresponds to a unique
orientation-reversing isometry taking $p$ to $y$. Conversely, this
isometry determines the plane. Horizontal planes contain none of the
lines. Every plane therefore contains at most $|S|\le n$ lines.
Apply the orientation-reversing case of
Theorem~\ref{tri:upd:thm:rich} to the proper sets $S$ and $US$,
both of size at most $n$, with overlap threshold $\sigma n$.
For such a plane, the number of its lines equals the number of points
of $S$ mapped into $US$ by the associated isometry. The cited theorem
therefore gives only $o_\sigma(n)$ planes containing at least
$\sigma n$ lines. Delete all lines in these planes. Since each plane
contains at most $n$ lines, this removes $o_\sigma(n^2)$ lines.

Each deleted line $\ell_{p,y}$ has a finite intersection with at most
$|S|-1$ other lines
of $\mathcal{L}_U$. Indeed, an intersecting line cannot have the same
source $p$ unless it is the same line, since the rotation represented
by an intersection has only one image of $p$. For a fixed source
$p'\ne p$, that rotation preserves the distance between the source
points, so an intersection requires $|p-p'|=|y-y'|$, and properness of $US$
allows at most one target $y'$. Each deleted line thus occurs in at
most $2(|S|-1)$ ordered intersecting pairs, counting its two possible
positions. The $o_\sigma(n^2)$ deletions remove only
$o_\sigma(n^3)$ such pairs.

We next show that a regulus contains at most $4M_2(S)$ lines of the
full family. Write a line as $(c,\tau)=(A+\tau B,\tau)$. Its direction and
moment vectors, with components $d_i$ and $m_i$, respectively, are
\[
 \mathbf{d}=(B_1,B_2,1),\qquad
 \mathbf{m}=(A_2,-A_1,A_1B_2-A_2B_1).
\]
Here $\mathbf{m}=(A_1,A_2,0)\times\mathbf{d}$, so these are the usual
\emph{Pl\"ucker coordinates} of the line. Since $p=A+JB$ and $y=A-JB$, the endpoint coordinates are
\[
 p=\left(\frac{-m_2-d_2}{d_3},
          \frac{m_1+d_1}{d_3}\right),\qquad
 y=\left(\frac{-m_2+d_2}{d_3},
          \frac{m_1-d_1}{d_3}\right).
\]
Thus, each endpoint is a linear projective projection in Pl\"ucker
coordinates. Each ruling of a regulus is a conic in these coordinates.
For the model quadric $x_0x_3-x_1x_2=0$, its line coordinates are
quadratic homogeneous functions of the ruling parameter $[u:v]$.
For example, one ruling consists of the lines spanned by
$(u,0,v,0)$ and $(0,u,0,v)$; their exterior products have coordinates
among $0,u^2,uv,v^2$, up to signs. The other ruling has the same
form after interchanging coordinates. Projective changes of coordinates act linearly on Pl\"ucker
coordinates, so the same conclusion holds for every regulus.

On a ruling containing members of $\mathcal{L}_U$, at least one endpoint
projection is nonconstant. Indeed, $d_3\ne0$ on a nonempty open
part of that ruling. If both endpoint projections were constant
there, the recovery formulas would give the same line throughout
this open part, which is impossible. Choose a nonconstant endpoint
projection. Its image is a curve of degree at most two, and each fiber
contains at most two ruling parameters. To check the fiber bound,
fix an endpoint and choose a projective line through it that does
not contain the image curve. Pulling back its linear equation gives
a nonzero homogeneous polynomial of degree at most two in $u,v$.
The fiber is contained in its zero set and therefore has size at
most two. Any base point of an endpoint
projection has $d_3=0$ and is absent from $\mathcal{L}_U$, where $d_3=1$.
Composition with the orthogonal map $U$ preserves polynomial degree,
so $M_2(US)=M_2(S)$. The image therefore contains at most $M_2(S)$
allowed endpoints, whether the chosen projection records $p$ or $y$. Each ruling therefore contains at most $2M_2(S)$ lines,
and the two rulings together contain at most $4M_2(S)$ lines.
By~\eqref{eq:conic}, the latter quantity is at most $4\sigma n$.
The plane bound after deletion is even smaller, namely $\sigma n$.
This proves all assertions of the lemma.
\end{proof}

For fixed $U$, an assigned equal-distance quadruple determines the
ordered pair $(\ell_{p_1,Uq_1},\ell_{p_2,Uq_2})$ of intersecting lines.
Distinct quadruples give distinct ordered pairs, since each line
determines its endpoint pair. Thus, the deletion estimate in
Lemma~\ref{lem:plane-regulus-bounds} discards only $o_\sigma(n^3)$
assigned quadruples for this $U$, and the same order of loss holds
after summing over the four choices of $U$.

\subsubsection{Counting intersections by angular localization}
The remaining equal-distance quadruples are represented by selected
intersection points associated with the isometries $Uf$.
We now count these points.
\par

The following estimate of Guth and Katz, in the form stated by
Cardinal, Payne, and Solomon~\cite[Theorem~3.4]{CPS}, bounds the
number of intersections once concentrations in planes and reguli
have been removed.

\begin{standard}\label{lit:rich-lines}
Let $\mathcal{L}$ be a finite set of distinct affine lines in $\R^3$.
If no plane or regulus contains more than $B$ members of $\mathcal{L}$,
then the number of points incident to at least two members is
$O(|\mathcal{L}|^{3/2}+B|\mathcal{L}|)$, with an absolute implicit constant.
\end{standard}

For points incident to at least $k$ lines, where $k\ge3$ is an
integer, no bound on the number of lines in a regulus is needed.
More precisely, Guth and Katz~\cite[Theorem~4.5]{GuthKatz} show that
if at most $B$ members of $\mathcal L$ lie in any plane, then the
number of such points in $\R^3$ is
\[
 O\!\left(\frac{|\mathcal L|^{3/2}}{k^2}
       +\frac{B|\mathcal L|}{k^3}+\frac{|\mathcal L|}{k}\right).
\]
An analogue over any field $\mathbb F$ of characteristic zero or
greater than two is due to Pham, Shen, Tran, and
Xue~\cite[Theorem~1.4]{PhamShenTranXue}. Under the same plane
hypothesis, their bound for lines in $\mathbb F^3$ is
\[
 O\!\left(\frac{|\mathcal L|^{3/2}+B|\mathcal L|}{k^2}
       +\frac{|\mathcal L|}{k}\right),
\]
with the additional assumption
$|\mathcal L|<(\operatorname{char}\mathbb F)^2/6$ in positive
characteristic. The implicit constants in both estimates are absolute.

\begingroup
For fixed $U$, each selected point represents $Uf$ with $\kappa_S(f)=2$,
and hence at most two equal-distance quadruples. Using only the bound
$n^2$ on the number of remaining lines and the plane and regulus bounds, however,
Theorem~\ref{lit:rich-lines} gives an $O(n^3)$ bound. To obtain the required saving,
we seek line families of size $o(n^2)$ with bounded overlap, such that the two lines
through every selected point belong together to at least one family.
An arbitrary subdivision of the lines need not preserve this property.

The other isometry associated with the same two matches provides the
required angular relation. Since $g$ maps the same source endpoints to
the same target endpoints and satisfies $|g(0)|=O_\theta(h)$, the two
endpoint pairs have nearly the same angular label. These labels will be
defined from the endpoints of each line, independently of the selected
isometry. Short overlapping intervals of labels keep each selected pair
of lines together, while angular nonconcentration makes each line family small.
We continue to count the points representing $Uf$; the role of $g$ is to
establish this relation between their two incident lines.
\par\endgroup

\begin{lemma}\label{lem:angular-localization}
After the deletions in Lemma~\ref{lem:plane-regulus-bounds}, the number of
equal-distance quadruples with a good source segment whose selected
realizing isometry $f$ satisfies $\kappa_S(f)=2$
and whose two associated
lines remain in the family for their chosen $U$ is at most
\[
 o_\theta(n^3)+C\sigma n^3,
\]
where $C$ is absolute.
\end{lemma}

\begin{proof}
Fix $U\in\{\mathrm{Id},-\mathrm{Id},R,-R\}$ and consider the
equal-distance quadruples assigned to
this choice, so that $Uf$ is a nonidentity rotation. Each selected $f$
has exactly two matched source points $p_1,p_2$, determined up to order;
write $q_i=f(p_i)$. Let $r$ be reflection in the line through $p_1,p_2$.
The isometry
\[
 g=f\circ r
\]
also takes $p_i$ to $q_i$ and has orientation opposite to that of $f$.
It is therefore the other isometry in Lemma~\ref{lem:two-isometries}. Since $f$ is the
previously selected isometry with $|f(0)|\ge c_\theta$, that lemma gives
$|g(0)|\le C_\theta h$ for all sufficiently large $j$.

Because $Uf$ preserves orientation, the sign
\[
 s:=\det O_g=-\det U
\]
is fixed for this choice of $U$. Assign to every line $\ell_{p,Uq}$
of the full family the angular label
\[
 b_s(p,q):=\arg q-s\arg p\quad\text{in }\R/(2\pi\mathbb Z).
\]
This label uses the original target $q$, before applying $U$, and depends
only on the endpoint pair and $U$. In particular, it does not depend on
which selected intersection the line is being considered at.

For the companion $g$ of a selected $f$, there is an angular parameter
$\beta_g$ such that
\[
 \arg(O_gx)=\beta_g+s\arg x\pmod{2\pi}\qquad(x\ne0).
\]
Writing $q_i=O_gp_i+g(0)$, and using the lower bound $|p_i|\ge1-h>1/2$,
we obtain
\begin{equation}\label{eq:companion-angular-label}
 b_s(p_i,q_i)=\beta_g+O_\theta(h)\pmod{2\pi},\qquad i=1,2.
\end{equation}
Here adding a vector of length $O_\theta(h)$ to a vector whose length is
bounded below changes its angle by $O_\theta(h)$. The parameter $\beta_g$
may vary with $f$, but the same parameter occurs for its two matches.
All angular error estimates are in circular distance.

Set $\omega=\sqrt{h+1/n}$. Then $\omega\to0$ and
$h/\omega\le\sqrt h\to0$. Partition the angle circle into
$\lceil2\pi/\omega\rceil$ equal half-open intervals. For sufficiently
large $j$, each interval has length between $\omega/2$ and $\omega$.
Enlarge each by its two neighbors on either side. With $\theta$ fixed,
\eqref{eq:companion-angular-label} implies that both labels of each selected
quadruple lie in the enlargement of the interval containing $\beta_g$.

For each enlarged interval $I$, let $\mathcal L_I$ be the remaining
lines with labels in $I$. Each line belongs to at most five such families,
and consequently
\[
 \sum_I|\mathcal L_I|\le5|S|^2\le5n^2.
\]
For fixed $p$, the condition $b_s(p,q)\in I$ is equivalent to
$\arg q\in s\arg p+I$. This is an interval of length at most $5\omega$,
so it contains at most $n\Omega_j(5\omega)$ possible targets in $S$.
Thus
\[
 \max_I|\mathcal L_I|\le n^2\Omega_j(5\omega)=o(n^2).
\]
The last limit follows from angular nonconcentration and $\omega\to0$.
These bounds are uniform in $S$ and in all subsequent subfamilies of
the remaining lines.

By Lemma~\ref{lem:plane-regulus-bounds}, each $\mathcal L_I$ inherits the
bound $4\sigma n$ on the number of lines in a plane or regulus. Apply Theorem~\ref{lit:rich-lines} to every nonempty
$\mathcal L_I$. Each selected intersection point is incident to exactly
two lines even in the full family $\mathcal L_U$, because
$\kappa_S(f)=2$. Ordering these two lines gives at most two equal-distance
quadruples. The incidence estimate may count additional points,
which only increases the upper bound. Hence the number of selected
quadruples for this $U$ is at most
\[
 \begin{aligned}
 C\sum_I\bigl(|\mathcal L_I|^{3/2}+4\sigma n|\mathcal L_I|\bigr)
 &\le C\sqrt{\max_I|\mathcal L_I|}\sum_I|\mathcal L_I|+C'\sigma n^3\\
 &\le C''n^3\sqrt{\Omega_j(5\omega)}+C'\sigma n^3\\
 &=o_\theta(n^3)+O(\sigma n^3).
 \end{aligned}
\]
The multiplicity bound uses $\kappa_S(f)=2$; no further estimate for
$\kappa_S(g)$ is used in this angular-localization argument.
The constants multiplying $\sigma$ are absolute. Dependence on $\theta$
enters through how large $j$ must be for the angular covering to apply.
Summing over the four choices of $U$ proves the lemma.
\end{proof}

\begin{proof}[Conclusion of the proof of Proposition~\ref{prop:annulus}]
Fix $\sigma>0$ and $0<\theta<\pi/8$. The initial restriction and
conic removals cost $o_\sigma(n^3)$. The discarded equal-distance
quadruples whose source segment is not good contribute
$O(\Omega_j(2\theta)n^3)$. Lemma~\ref{lem:far} bounds the contribution
of quadruples whose selected realizer $f$ satisfies
$\kappa_S(f)\ge3$ by $o_\theta(n^3)$, while
Lemmas~\ref{lem:plane-regulus-bounds}
and~\ref{lem:angular-localization} handle the remaining quadruples.
Together, these estimates give

\[
 \frac{E(P_j)}{|P_j|^3}
 \le o_\sigma(1)+o_\theta(1)
      +C\Omega_j(2\theta)+C'\sigma,
\]
where $C,C'$ are absolute. Letting $j\to\infty$, with $\sigma$ and
$\theta$ fixed, yields
\[
 \limsup_{j\to\infty}\frac{E(P_j)}{|P_j|^3}
 \le C\limsup_{j\to\infty}\Omega_j(2\theta)+C'\sigma.
\]
Next let $\sigma\downarrow0$, and then let $\theta\downarrow0$.
The angular nonconcentration condition~\eqref{eq:angular-nonconcentration}
gives $E(P_j)=o(|P_j|^3)$. Finally, $|P_j|\ge2$ and $E(P_j)>0$ for
all sufficiently large $j$, so division in \eqref{eq:cs-distances}
is legitimate and gives
\[
 \frac{D(P_j)}{|P_j|}
 \ge \frac{(1-1/|P_j|)^2}{E(P_j)/|P_j|^3}\longrightarrow\infty.
\]
This completes the proof of Proposition~\ref{prop:annulus}.
\end{proof}

\section{Sharpness examples of Proposition \ref{prop:annulus}}\label{sec:annulus-examples}

The following two lemmas show that both properness and angular
nonconcentration are necessary for the subcubic energy estimate in
Proposition~\ref{prop:annulus}. The first example also shows that
angular nonconcentration alone does not imply a superlinear number
of distinct distances.

\begin{lemma}\label{lem:annulus-without-properness}
For each integer $j\ge1$, let $n_j=2j+1$, and let $P_j$ be the
vertices of a regular $n_j$-gon on $S^1$. These sets satisfy the
radial and angular hypotheses of Proposition~\ref{prop:annulus}
with $P'_j=P_j$ and $h_j=\frac{1}{n_j}$, but they are not proper. Moreover,
\[
 D(P_j)=\frac{n_j-1}{2},\qquad
 E(P_j)=2n_j^2(n_j-1).
\]
In particular,
\[
 \frac{D(P_j)}{|P_j|}\longrightarrow\frac12,\qquad
 \frac{E(P_j)}{|P_j|^3}\longrightarrow2.
\]
\end{lemma}

\begin{proof}
All points lie on $S^1$, so the radial hypothesis holds with the
stated choices of $P'_j$ and $h_j$. The angles of consecutive
vertices differ by $2\pi/n_j$. Thus, every angular interval $I$
satisfies
\[
 \bigl|\{x\in P_j:\arg x\in I\}\bigr|
 \le\frac{n_j|I|}{2\pi}+1.
\]
Given $\varepsilon>0$, choose $0<t\le2\pi$ with
$t/(2\pi)\le\varepsilon/2$. For all sufficiently large $j$,
we also have $1/n_j\le\varepsilon/2$. The displayed bound then
gives the required angular nonconcentration for every interval
with $|I|\le t$.

At each vertex, the two adjacent vertices are equidistant from
it, so $P_j$ is not proper. Its distinct positive squared distances
are
\[
 4\sin^2\frac{\pi k}{n_j},\qquad
 1\le k\le\frac{n_j-1}{2}.
\]
These values are distinct because sine is strictly increasing
on $(0,\pi/2)$. For each such value, every vertex has exactly two
vertices at that squared distance. Its ordered multiplicity is
therefore $2n_j$. Consequently,
\[
 D(P_j)=\frac{n_j-1}{2},\qquad
 E(P_j)=\frac{n_j-1}{2}(2n_j)^2=2n_j^2(n_j-1).
\]
Dividing by $|P_j|=n_j$ and $|P_j|^3=n_j^3$, respectively, gives
the asserted limits. This completes the proof.
\end{proof}

\begin{lemma}\label{lem:annulus-without-nonconcentration}
For each integer $m\ge2$, set
\[
 \begin{gathered}
 r_i=\sqrt{1+\frac{i}{m^2}}\quad(1\le i\le m),\\
 P_m=\{(r_i,0):1\le i\le m\}
       \cup\{(0,r_i):1\le i\le m\},\qquad n=|P_m|=2m.
 \end{gathered}
\]
The sets $P_m$ are proper and satisfy the radial hypothesis of
Proposition~\ref{prop:annulus} with $P'_m=P_m$ and $h_m=\frac{1}{2m}$,
but they do not satisfy angular nonconcentration. Moreover,
\[
 E(P_m)>\frac{n^3}{4},\qquad
 D(P_m)=\binom{m}{2}+2m-1.
\]
Thus, the subcubic energy conclusion fails, although
$D(P_m)/|P_m|\to\infty$.
\end{lemma}

\begin{proof}
For every $1\le i\le m$,
\[
 1<r_i\le\sqrt{1+\frac1m}<1+\frac1{2m}.
\]
The radial hypothesis therefore holds with the stated choices.
Exactly $m=n/2$ points have angle $0$, and the other $m$ have
angle $\pi/2$. For any $t>0$, an angular interval of length at
most $t$ containing $0$ contains at least $n/2$ points. Hence the
angular hypothesis fails, for example with $\varepsilon=1/4$.

We next prove that all positive distances between points on the
same ray are distinct. For $1\le i<i+k\le m$,
\[
 2m^2(r_{i+k}-r_i)=\frac{2k}{r_{i+k}+r_i}.
\]
Since $2<r_{i+k}+r_i<2+1/m$ and $k<m$, we have
\[
 k-\frac12
 < k-\frac{k}{2m+1}
 =\frac{2k}{2+1/m}
 <2m^2(r_{i+k}-r_i)<k.
\]
The intervals $(k-1/2,k)$ are disjoint for distinct positive
integers $k$, so a distance determines $k$. For fixed $k$, the
quantity $r_{i+k}-r_i$ strictly decreases with $i$, because its
denominator $r_{i+k}+r_i$ strictly increases. Thus, each ray
determines exactly $\binom{m}{2}$ distinct positive distances,
and the two rays determine the same set of such distances.

Every distance within a ray is smaller than $1/(2m)$, whereas
every distance between the rays is greater than $\sqrt2$.
For a fixed point $(r_i,0)$, the squared distances to the other
ray are
\[
 |(r_i,0)-(0,r_k)|^2=2+\frac{i+k}{m^2},\qquad 1\le k\le m,
\]
which are distinct as $k$ varies. The same argument applies with
the two rays interchanged. Together with the preceding bounds,
this shows that all distances from any fixed point are distinct.
Hence $P_m$ is proper.

There are $2m^2$ ordered pairs with endpoints on different rays.
Their squared distances take exactly $2m-1$ values, corresponding
to $i+k=2,\ldots,2m$. By the Cauchy--Schwarz inequality, these pairs alone give
\[
 E(P_m)\ge\frac{(2m^2)^2}{2m-1}
 >2m^3=\frac{n^3}{4}.
\]
The distances within the rays are disjoint from those between
the rays. Consequently,
\[
 D(P_m)=\binom{m}{2}+2m-1,\qquad
 \frac{D(P_m)}{|P_m|}
 =\frac{m-1}{4}+1-\frac1{2m}\longrightarrow\infty.
\]
This proves all assertions of the lemma.
\end{proof}

Lemma~\ref{lem:annulus-without-nonconcentration} concerns the
energy conclusion of Proposition~\ref{prop:annulus}. A proper
sequence for which $|P_m|\to\infty$ and $D(P_m)/|P_m|$ remains
bounded would instead answer the Erd\H{o}s--Davies divergence
question in the negative. Conversely, any such sequence could
be placed in disks of radii tending to zero about $(1,0)$ by
similarities, which preserve properness, $D(P_m)$, and $E(P_m)$.
\par

\section{An open question on polynomial extensions}
\label{sec:polynomial-conjecture}

In this section, we propose an open question for symmetric
polynomials.

Let $F:\R^2\times\R^2\to\R$ be a fixed nonconstant symmetric
polynomial, so that $F(p,q)=F(q,p)$ for all $p,q\in\R^2$.
We say that a finite set $P\subset\R^2$ is \emph{admissible
for $F$} if $q\mapsto F(p,q)$ is injective on $P$ for each $p\in P$,
and we write

\[
 V_F(P)=\{F(p,q):p,q\in P\}.
\]
For $F(p,q)=|p-q|^2$, admissibility is precisely properness, and
$|V_F(P)|=D(P)+1$ for nonempty $P$.

Our question is stated as follows.

{\bf Question.}
Which such polynomials $F$ admit a constant $K>0$ and admissible sets
$P_j\subset\R^2$ such that
\[
 |P_j|\longrightarrow\infty,
 \qquad |V_F(P_j)|\le K|P_j|?
\]

The following proposition gives a sufficient condition for
the existence of such sequences.

\begin{proposition}\label{prop:poly-sum}
Let $I$ be an open interval, and let $\gamma:I\to\R^2$ be an injective
real-analytic map with $\gamma'\ne0$. Suppose that there is a nonconstant
real-analytic function $h:I+I\to\R$ such that
\begin{equation}\label{eq:poly-sum}
 F(\gamma(s),\gamma(t))=h(s+t)\qquad(s,t\in I).
\end{equation}
Then, for every $n\ge1$, there is an admissible $n$-point set $P$
with $|V_F(P)|=2n-1$.
\end{proposition}
\begin{proof}
Since $I+I=2I$ and $h$ is nonconstant, there is $t_0\in I$
with $h'(2t_0)\ne0$. Choose an open interval
$J\subset I$ about $t_0$ such that $h$ is injective on $J+J$.
Take an $n$-term arithmetic progression $A\subset J$ and put
$P=\gamma(A)$. Equation~\eqref{eq:poly-sum} and the injectivity of
$h$ give admissibility. Moreover,
$|V_F(P)|=|A+A|=2n-1$, as required.
\end{proof}

For $F(p,q)=p_1+q_1$, identity~\eqref{eq:poly-sum} holds with
$\gamma(t)=(t,0)$ and $h(u)=u$.
For the dot product $F(p,q)=p\cdot q$, it holds with
$\gamma(t)=(2^t,0)$ and $h(u)=2^u$, corresponding to the
geometric-progression construction in~\cite[Example~1.2]{ValculescuDeZeeuw}.

It would be very interesting to know whether the sufficient condition in Proposition~\ref{prop:poly-sum}
also necessary for the existence of the sequences in the above question.

Note that if $F(p,p)$ is constant, then \eqref{eq:poly-sum} would force $h$
to be constant on $I+I$, since
$h(2t)=F(\gamma(t),\gamma(t))$ for $t\in I$.
Thus, for squared distance, an affirmative answer to this
converse question would imply the Erd\H{o}s--Davies divergence
conjecture.

\end{document}